\documentclass[11 pt]{article}
\usepackage[top=2cm, bottom=2cm, left=25mm, right=25mm]{geometry}

\usepackage[utf8]{inputenc}
\usepackage{amsmath}
\usepackage{parskip}
\usepackage{amssymb}
\usepackage{color}
\usepackage[all]{xy}  
\usepackage{tikz-cd}
\usepackage{float}
\usepackage{setspace}
\usepackage{hyperref}
\newtheorem{theorem}{Theorem}

\newtheorem{proposition}[theorem]{Proposition}

\newtheorem{definition}{Definition}

\usepackage{caption}
\usepackage{lineno}

\usepackage{algorithmicx}
\usepackage{algorithm}
\usepackage{algpseudocode}

\usepackage[acronym]{glossaries}
\newacronym{RKMK}{RKMK}{Método de Runge-Kutta-Munthe-Kaas}

\usepackage{hyperref}
\hypersetup{colorlinks = true}

\usepackage{tikz-cd}
\usepackage{pgfplots}
\usepgfplotslibrary{groupplots}
\pgfplotsset{compat = 1.16}
\usetikzlibrary{decorations.markings}

\newcommand{\Sk}{{\rm\ \!S}}            
\newcommand{\Ck}{{\rm\ \!C}}           
\newcommand{\Tk}{{\rm\ \!T}}             
\newcommand{\Vk}{{\rm\ \!V}}             

\newcommand{\kk}{\kappa}

\newcommand{\bfI}{{\mathbf I}_{\kap}}

\def\k{{\kappa}}

\newcommand{\te}{\phi}

\numberwithin{equation}{section}
\numberwithin{theorem}{section}
\numberwithin{definition}{section}
\def\k{{\kappa}}
 \def \kap{{\boldsymbol{\kappa}}}

\author{L. Blanco$^{*}$, F. Jiménez$^{*}$, J. de Lucas$^\dagger$, C. Sardón$^{*}$}
\begin{document}
\centerline{\Large Geometry preserving numerical methods}
\vskip .3cm
\centerline{\Large for physical systems with finite-dimensional Lie algebras}
\vskip .5cm

\centerline{L. Blanco$^{*}$, F. Jiménez$^{\S}$, J. de Lucas$^\dagger$, C. Sardón$^{*}$}
\vskip .2cm

\begin{center}
 $^\dagger$Centre de Recherches Math\'ematiques, 
 Universit\'e de Montr\'eal, \\
 \medskip 
 Pavillon André-Aisenstadt, Chemin de la Tour 2920,\\
 \medskip
Montréal (Québec) H3T 1J4, Canada.\\
\medskip 
$^{\dagger}$Department of Mathematical Methods in Physics, University of Warsaw,\\
\medskip
	ul. Pasteura 5, 02-093, Warsaw, Poland.\\
 \medskip
\end{center}
	\medskip
	\centerline{$^{*}\!\!$ Departamento de Matem\'aticas Aplicadas a la Ingenier\'ia Industrial,}
 \medskip
\centerline{Universidad Polit\'ecnica de Madrid (UPM)}
	\medskip
	\centerline{c. José Gutiérrez Abascal 2, 28006, Madrid, Spain.}
\medskip
\centerline{$^{\S}\!\!$ Departamento de Matem\'atica Aplicada I, ETSII,}
 \medskip
\centerline{Universidad Nacional de Educaci\'on a Distancia (UNED)}
	\medskip
	\centerline{c. Juan del Rosal 12, 28040, Madrid, Spain.}
	\medskip

\begin{abstract}
We propose a geometric integrator to numerically approximate the flow of Lie systems. The key is a novel procedure that integrates the Lie system on a Lie group intrinsically associated with a Lie system on a general manifold via a Lie group action, and then generates the discrete solution of the Lie system on the manifold via a solution of the Lie system on the Lie group.

One major result from the integration of a Lie system on a Lie group is that one is able to solve all associated Lie systems on manifolds at the same time, and that Lie systems on Lie groups can be described through first-order systems of linear homogeneous ordinary differential equations (ODEs) in normal form. This brings a lot of advantages, since solving a linear system of ODEs involves less numerical cost. Specifically, we use two families of numerical schemes on the Lie group, which are designed to preserve its geometrical structure: the first one based on the Magnus expansion, whereas the second is based on Runge-Kutta-Munthe-Kaas (RKMK) methods. Moreover, since the aforementioned action relates the Lie group and the manifold where the Lie system evolves, the resulting integrator preserves any geometric structure of the latter. We compare both methods for Lie systems with geometric invariants, particularly a class on Lie systems on curved spaces. We also illustrate the superiority of our method for describing long-term behavior and for differential equations admitting solutions whose geometric features depends heavily on initial conditions.

As already mentioned, our milestone  is to show that the method we propose preserves all the geometric invariants very faithfully, in comparison with nongeometric numerical methods.

\end{abstract}

{\it MSC 2020 classes: 34A26; 53A70 (primary) 37M15; 49M25 (secondary)}

\section{Introduction}

The history of numerical methods on Lie groups is intertwined with the development of computational mathematics and the study of Lie theory. 
The foundations of Lie theory were settled by the Norwegian mathematician Sophus Lie in the late 19th century. However, it was not until the 20th century that the application of Lie groups to practical problems and the development of numerical methods gained momentum.

 In the 1970s, mathematicians and physicists began to explore numerical integration methods for Lie group equations of motion. Afterwards, pioneering work by Blanes, Casas, Oteo, and Ros provided explicit symplectic integrators for specific Lie groups, such as the rotation group $SO(3)$ and the special Euclidean group $SE(3)$. These methods preserved important geometric properties of Lie groups, such as the energy and/or the symplecticity \cite{blanes06,blanes00}.

The computation of geodesics on Lie groups became a topic of interest in the 1980s. Researchers like Murray, Arimoto, and Sastry developed numerical methods to compute geodesics on Lie groups such as $SO(3)$ and $SE(3)$ \cite{Murray1,Murray2}. These methods relied on various techniques, including the exponential map, interpolation, and numerical optimization algorithms. The optimization of functions defined on Lie groups gained prominence in the 1990s. Researchers such as Absil, Mahony, and Mallick developed numerical optimization algorithms specifically tailored to the geometric properties of Lie groups \cite{Mahony1,Mahony2,Mahony3}. These methods allowed for efficient optimization of functions over Lie groups, which found applications in robotics, computer vision, and control theory.
The interpolation of motions on Lie groups received significant attention in the early 2000s. Researchers like Sola, Kuffner, and Agrawal proposed interpolation algorithms for Lie group elements, enabling smooth and visually appealing motion planning in applications such as robotics and computer graphics \cite{Sola}.

 In recent years, there has been continued progress in numerical methods on Lie groups, fueled by advancements in computational power and the increasing demand for efficient algorithms in applications. Research continues to focus on refining existing methods, developing new techniques, and exploring applications in areas like machine learning, motion planning, and optimization.
 
The Runge-Kutta methods are a family of numerical integration techniques commonly used to solve ordinary differential equations. They involve evaluating the derivative of the function at multiple points within a time step and using a weighted sum of these derivatives to update the solution. A comprehensive survey on modern geometric Lie group methods, including new ideas and techniques, can be found in \cite{iserles05}.

The RKMK method combines these two concepts by using the Munthe-Kaas rule to select the sampling points in the Runge-Kutta integration scheme. RKMK method is also the term we use to refer to the usual Runge-Kutta method (RK) applied on Lie groups. By considering the distribution of the highest derivative of the function being integrated, the RKMK method aims to improve the accuracy and efficiency of the integration process \cite{munthe98,munthe99,munthe96}. 

The specific details of the RKMK method, including the choice of sampling points and the weights assigned to the derivatives, can vary depending on the implementation and the problem at hand. Researchers have proposed different variants of the RKMK method with varying degrees of accuracy and computational complexity. Since the properties of a RKMK methods are the same of a classical RK, the symplecticity is preserved for certain orders: for example, the second-order Störmer-Verlet method, also known as the leapfrog method, is a well-known second-order symplectic integrator \cite{Verlet}. There are several fourth-order symplectic integrators, such as the Forest-Ruth method and the Yoshida method \cite{Yoshida}.
Higher-order symplectic integrators have also been developed, such as the sixth-order McLachlan integrator \cite{McLa} and the eighth-order Blanes-Moan integrator \cite{Moan}.
These symplectic Runge-Kutta methods are designed to preserve the symplectic structure of Hamiltonian systems and offer improved accuracy and long-term stability compared to non-symplectic methods.

It is important to note that the choice of a specific symplectic Runge-Kutta method depends on the requirements of the problem at hand, including the desired accuracy, computational efficiency, and preservation of particular properties. In our case, we will work with a fourth-order RKMK.

\textbf{The RKMK Methods}

The basic idea behind applying the fourth-order RKMK method is to update the group elements using Lie group operations while approximating the derivatives of the group elements at multiple intermediate points within a time step. The following steps outline a typical approach.
\begin{itemize}
\item[--] Initialization: Start with an initial group element.

\item[--]  Time Step Selection: Choose an appropriate time step size for the integration process.

\item[--]  Derivative Evaluation: Evaluate the derivative of the group element at the initial time.

\item[--]  State Update: Use the fourth-order RK method to update the group element by integrating the derivative. This involves evaluating the derivative at multiple intermediate points within the time step and combining them with weighted sums to update the state.

\item[--]  Group Operation: Apply appropriate Lie group operations (e.g., matrix multiplication, exponentiation) to ensure the updated state remains on the Lie group manifold.

\item[--]  Repeat: Repeat steps 3-5 until the desired integration time is reached.

\end{itemize}
By incorporating the Lie group operations in the state update step and properly handling the derivatives, the fourth-order RK method can be applied to approximate solutions on Lie groups.

\textbf{Magnus Method and its Interpretation}

To solve the initial-value problem for the linear system of ordinary differential equations on $\mathbb{R}^n$ of the form
\begin{equation}\label{eq:InValProb}
\frac{dY}{dt}(t) = A(t)Y(t), \quad Y(t_0) = Y_0,
\end{equation}
where \(Y(t)\) is an unknown \(n\)-dimensional vector $t$-dependent function and \(A(t)\) is an \(n \times n\) matrix with $t$-dependent entries, the Magnus approach was introduced. The solution for \(n = 1\) is 
\[Y(t) = \exp \left(\int_{t_0}^t A(s) \, ds\right) Y_0.\]
This solution also holds for \(n > 1\) provided \(A(t_1)A(t_2) = A(t_2)A(t_1)\) for any pair of values \(t_1\) and \(t_2\), especially when \(A\) is independent of \(t\). However, for the general case, the aforementioned expression is not a valid solution.

Wilhelm Magnus devised a method to solve the matrix initial-value problem (\ref{eq:InValProb}) by introducing the exponential of a specific \(n \times n\) matrix function \(\Omega(t, t_0)\) as follows
\begin{equation}\label{expeq}
Y(t) = \exp \left(\Omega(t, t_0)\right) Y_0,
\end{equation}
where \(\Omega(t)\) is constructed as a series expansion
\[\Omega(t) = \sum_{k=1}^\infty \Omega_k(t),\]
with \(\Omega(t)\) representing \(\Omega(t, t_0)\) for simplicity and taking \(t_0 = 0\).

Using \(\frac{d}{dt}(e^\Omega) e^{-\Omega} = A(t)\) and the Poincaré-Hausdorff matrix identity, Magnus related the time derivative of \(\Omega\) to the generating function of Bernoulli numbers and the adjoint endomorphism of \(\Omega\), namely ${\rm ad}(\Omega)=[\Omega,\cdot]$, as follows
\begin{equation}\label{eq:MagExp}
\frac{d\Omega}{dt} = \frac{\operatorname{ad}(\Omega)}{\exp(\operatorname{ad}(\Omega))-{\rm Id}_n} A
\end{equation}
to solve \(\Omega\) in a recursive manner in terms of \(A\) as a continuous analog of the Baker-Campbell-Hausdorff (BCH) expansion \cite{Hall2,varadajan}. Note that ${
\rm Id}_n$ stands for the $n\times n$ identity matrix.

Equation (\ref{eq:MagExp}) is called the {\it Magnus expansion}, or {\it Magnus series}, for the solution of (\ref{eq:InValProb}). The first four terms of this series read
\begin{align*}
\Omega_1(t) &= \int_0^t A(t_1) \, dt_1, \\
\Omega_2(t) &= \frac{1}{2} \int_0^t dt_1 \int_0^{t_1} [A(t_1), A(t_2)]dt_2 , \\
\Omega_3(t) &= \frac{1}{6} \int_0^t dt_1 \int_0^{t_1} dt_2 \int_0^{t_2}  \left([A(t_1), [A(t_2), A(t_3)]] + [A(t_3), [A(t_2), A(t_1)]]\right)dt_3, \\
\Omega_4(t) &= \frac{1}{12} \int_0^t dt_1 \int_0^{t_1} dt_2 \int_0^{t_2} dt_3 \int_0^{t_3}  \left([[[A(t_1), A(t_2)], A(t_3)], A(t_4)] + \ldots \right)dt_4.
\end{align*}

By expressing the solution in terms of the exponential of a matrix function \eqref{expeq}, the Magnus series offers a systematic way to approximate the solution. The Magnus approach very often preserves important qualitative properties of the exact solution, such as the symplectic or unitary character, even in truncated forms. This method has found applications in various fields, including classical mechanics and quantum mechanics, where it offers an alternative to conventional perturbation theories. The Magnus expansion method stands as a valuable tool for analyzing and approximating solutions to linear differential equations \cite{blanes00,Magnus}, and, naturally, has strong applications when $Y(t)$ belongs to a matrix Lie group.

These two numerical integration techniques play a crucial role in the theory of Lie systems. These methods aim to preserve the geometric structures and qualitative properties of the underlying system, such as symplecticity, conservation laws and, in the first place, the Lie group structure itself. The Magnus expansion and RKMK methods are particularly useful for preserving the long-term behavior of Hamiltonian systems, and, in particular, the so called Lie-Hamilton systems \cite{araujo20}, which are a special class of Lie systems. 

Hamiltonian systems relative to a symplectic structure admit natural constants of motion, and numerical methods that preserve such symplectic structures and related invariants ensure that obtained approximate solutions will remain  at every moment in certain regions containing the searched exact solution, e.g., in submanifolds given by a constant of motion. In the short term, this may not have a big impact in the accuracy of the chosen numerical solutions. But in the long term, this may ensure that the obtained numerical solution will always remain in a specific region containing the exact solution. In particular, symplectic preserving numerical methods may ensure that a numerical solution around a stable point of a Hamiltonian system given by the minimum of its Hamiltonian function remains close to the equilibrium point. Other numerical methods may not be able to reproduce such a behavior.

Lie systems occurred for the first time in the study of Riccati equations \cite{Ko83} as a consequence of the generalisation to a nonlinear realm of the known superposition rules for nonautonomous linear systems of first-order ordinary differential equations. Among other reasons, superposition rules are interesting to solve numerically systems of differential equations whose general solutions cannot be exactly found \cite{Wi_83}. Although most differential equations cannot be studied via Lie systems, Lie systems have many relevant applications in physics, control theory, and other fields \cite{lucas11,araujo20}. In particular, Lie-Hamilton systems occur in the study of Smorodinsky-Winternitz oscillators, Milne--Pinney equations, dissipative harmonic oscillators,  trigonometric oscillators, and so on (see \cite{araujo20} and references therein). Certain quantum mechanical systems, like quantum mechanical oscillators with time-dependent frequency and other time-dependent parameters, can also be studied via Lie systems on Lie groups \cite{lucas11}. Particular cases of matrix Riccati equations, which are also Lie systems, are associated with Painlev\'e trascendents, Sawada-Kotera equations, Kaup–Kupershmidt equations, etcetera \cite{LG17}. For all that reasons, the study of Lie systems is fully motivated from the point of view of applications.

Our approach to the numerical integration of Lie systems in this manuscript is the Lie group approach, which is relevant since there exists a Lie group action establishing a relationship between a certain Lie group intrinsically associated to the Lie system and the manifold where the Lie system itself evolves \cite{carinena00}. The two Lie group integrators that we have introduced exploit the algebraic structure of the Lie group associated with the Lie system to construct accurate and efficient numerical schemes. 

Our aim in this work is to depart from this preexisting technology on Lie group integrators and take advantage of it when numerically integrating Lie systems. These can evolve on manifolds with an additional particular compatible geometric structure, such as a group structure or curvature, and, therefore, a geometric integrator for them is in order. The Lie group action relating the Lie group underlying a Lie system and the manifold where it evolves represents a perfect tool to achieve this goal, generating a discrete sequence of points in the manifold (which naturally inherit its geometry) from a discrete sequence of points in the Lie group, which can be obtained from the Lie group integrators. In this way we establish a novel geometric integrator \cite{hairer06}, which we will test on a particular class of Lie systems on curved spaces and other Lie systems with interesting properties. Our focus in this article, as shall be noticed during its reading, is the dynamical (convergence) and geometrical  behavior of the discrete solutions of our integrators. Therefore, we have not paid much attention to the computation cost and time, which is a relevant issue to address in the future. However, as in other geometric integrators \cite{hairer06}, it is to be expected a higher computation cost, due to the their particular design accomodating the intrinsic nature of the systems they approximate. On the other hand  this special design shall redound on the long-term behavior, as pointed out in the next paragraph and throughout the paper.  Note also that ordinary differential equations in normal form around equilibrium points may be approximated by linear systems of differential equations, which can be studied via Lie systems \cite{araujo20}. It may happen that higher-order approximations may be described through Lie systems. This also justifies the potential applicability of our techniques.

The previous approach has natural advantages to describe the long-term behavior of obtained numerical solutions. Indeed, Lie systems have an evolution that is determined by a family of vector fields that span an integrable generalized distribution. Then, the Lie system has a related stratification of the manifold and the solutions of the Lie system always evolve within them \cite{carinena00}. By using numerical methods in the Lie group associated with the Lie system, the numerical solutions obtained will always remain within the stratification related to the Lie system. In a natural manner, this ensures that obtained numerical solutions will always be well behaved not escaping away from the strata containing the exact solution to be obtained. In general, this may not have a deep impact in the short-term behavior. Nevertheless, it may be relevant in solving problems with a strong dependence on initial conditions, as shown in the examples of this work, which have some kind of ``chaotic/bifurcation" behavior. Moreover, our methods ensure that numerical solutions will always satisfy certain geometric properties of the searched exact solutions  in the long-term. In contrast, general numerical methods may provide numerical solutions whose features in the long-term have nothing to do with the searched solution. 

The outline of the paper goes as follows. In \S\ref{GeomFun} we introduce the fundamentals on Lie groups and Lie algebras needed hereafter. Moreover, we describe  automorphic Lie systems and how to solve them in its underlying Lie group. The definition of the action of a Lie group on the manifold where the Lie system evolves is also presented. These two elements allow for the definition of the 7-step method for the reduction procedure of automorphic Lie systems in Definition \ref{Metodo7Pasos}. \S\ref{cap.discre} depicts some basics on numerical schemes and the Lie group methods employed afterwards. In \S\ref{metodos} we combine all the previous elements to propose our geometric method to  numerically integrate  automorphic Lie systems in Definition \ref{MetodosLS}. In \S\ref{EjemploDetallado} we pick a class of Lie systems on curved spaces, we apply with high detail the 7-step method to analytically solve them and, afterwards, we employ our integrator, showing its geometric properties. The advantages of our methods in the description of problems with a strong dependence on the initial conditions and with respect to long-term analysis of solutions is analyzed in Section \S\ref{Sec:IniLongNum}. The conclusions and outlook of our work are described in Section \S\ref{Sec:ConOut}.

\section{Geometric fundamentals}\label{GeomFun}

\subsection{Lie groups and matrix Lie groups} \label{fundaLie}

Let $G$ be a Lie group and let $e$ be its neutral element. Every  $ g \in G $ defines a right-translation $R_g : h\in G \mapsto hg\in G $ and a left-translation $L_g:h\in G\mapsto gh\in G$ on $G$.
A vector field, $ X^\mathrm{R} $, on $ G $ is {right-invariant} if $ X^\mathrm{R}(hg) = R_{g*,h} X^\mathrm{R}(h) $ for every $ h, g \in G $, where $ R_{g*, h} $ is the tangent map to  $ R_g $ at $h\in G$. The value of a right-invariant vector field, $ X^\mathrm{R} $, at every point of $ G $ is determined by its value at $e$, since, by definition, $ X^\mathrm{R}(g) = R_{g*,e} X^\mathrm{R}(e) $ for every $ g \in G$. Hence, each right-invariant vector field $X^R$ on $ G $ gives rise to a unique $X^R(e)\in T_eG$ and vice versa. Then, the space of right-invariant vector fields on $G$ is finite-dimensional, and it can be proved to be also a Lie algebra.
Similarly, one may define {left-invariant} vector fields on $G$, establish a Lie algebra structure on the space of left-invariant vector fields and set an isomorphism between the space $\mathfrak{g}$ of left-invariant vector fields on $G$ and $T_eG$. The Lie algebra of left-invariant vector fields on $G$, with the  Lie bracket $ [\cdot \, , \cdot] : \mathfrak{g} \times \mathfrak{g} \to \mathfrak{g} $ given by the commutator of vector fields, induces in $T_eG$ a Lie algebra via the identification of left-invariant vector fields with their values at $e$. Note that we will frequently identify $\mathfrak{g}$ with $T_eG$ to simplify the notation.

There is a natural mapping from $\mathfrak{g}$ to $G$, the so-called {exponential map}, of the form $ \exp :a\in  \mathfrak{g} \mapsto \gamma_a(1)\in G $, where $ \gamma_a: \mathbb{R} \to G $ is the integral curve of the right-invariant vector field $ X^\mathrm{R}_a $ on $G$ satisfying $ X^\mathrm{R}_a(e) = a $ and $ \gamma(0) = e $. It is worth noting that the exponential map could be defined using left-invariant vector fields to give exactly the same mapping (cf. \cite{Hall2}). If $\mathfrak{g}=\mathfrak{gl}(n,\mathbb{K})$, where $\mathfrak{gl}(n,\mathbb{K})$ is the Lie algebra of $n\times n$ square matrices with entries in a field $\mathbb{K}$ relative to the Lie bracket given by the commutator of matrices, then $\mathfrak{gl}(n,\mathbb{K})$ can be considered as the Lie algebra of the Lie group ${\rm GL}(n,\mathbb{K})$ of $n\times n$ invertible matrices with entries in $\mathbb{K}$. It can be proved that in this case $\exp:A\in \mathfrak{gl}(n,\mathbb{K})\mapsto \exp(A)\in {\rm GL}(n,\mathbb{K})$ is given by the standard expression of the exponential of a matrix \cite{lee00}, namely
\[ \exp(A) = {\rm I}_n + A + \frac{A^2}{2} + \frac{A^3}{6} + \cdots = \sum_{k = 0}^{\infty} \frac{A^k}{k!},\]
where we recall that ${\rm I}_n$ stands for the $n\times n$ identity matrix.

From the definition of the exponential map $\exp:T_eG\rightarrow G$, it follows that $\exp(sa)=\gamma_a(s)$ for each $ s \in \mathbb{R} $ and $a\in T_eG$. Let us show this. Given the right-invariant vector field $X^\mathrm{R}_{sa} $, with $ sa \in T_eG $, one has 
$$
X^\mathrm{R}_{sa}(g) = R_{g*,e} X^\mathrm{R}_{sa}(e) = R_{g*,e}(sa) = s R_{g*,e}(a),\qquad \forall g\in G.
$$
In particular, for $s=1$, it follows that $ X^\mathrm{R}_a(g) =R_{g*,e}(a) $ and, for a general $s$, one has that  $X^\mathrm{R}_{sa} = s X^\mathrm{R}_a $. Hence, if $ \gamma_a,\gamma_{sa}:\mathbb{R}\rightarrow G$ are the integral curves of $ X^\mathrm{R}_a $ and $X^\mathrm{R}_{sa}$ with initial condition $ e $, then it can be proved that, for $u=ts$, one has that
$$
\frac{d}{dt}\gamma_a(ts)=s\frac{d}{du}\gamma_a(u)=sX^{\rm R}_a(\gamma_a(ts))=X^{\rm R}_{sa}(\gamma_a(ts)).
$$
and $t\mapsto \gamma_a(st)$ is the integral curve of $X_{sa}^{\rm R}$ with initial condition $e$. Hence, $\gamma_{a}(st)=\gamma_{sa}(t)$.
Therefore, $ \exp(sa) = \gamma_{sa}(1) = \gamma_a(s)$. It is worth stressing that it is a consequence of
Ado's theorem \cite{ado47} that every Lie group can be written as a matrix Lie group on some open neighborhood of its neutral element. %If $ G $ is a matrix Lie group and $X\in \mathfrak{g} $, then \cite{lee00}

The exponential map establishes a diffeomorphism from an open neighborhood $U_\mathfrak{g}$ of $0$ in $T_eG$  and $\exp(U_\mathfrak{g})$. More in detail, every basis $ \mathcal{V} = \{ v_1, \dots,v_r \}$ of $T_eG$ gives rise to the so-called  {canonical coordinates of the second-kind} related to $\mathcal{V}$ defined by the local diffeomorphism
\begin{equation*}
\begin{array}{ccc}
U_\mathfrak{g}\subset T_eG                & \longrightarrow    & \exp(U_\mathfrak{g}) \subset G \\
(\lambda_1, \dots, \lambda_r) & \mapsto & \prod_{\alpha = 1}^{r} \exp( \lambda_\alpha v_\alpha ) \, ,
\end{array}
\end{equation*}
for an appropriate open neighborhood $U_\mathfrak{g}$  of $0$ in $T_eG\simeq\mathfrak{g}$. %a neighborhood of the neutral element in $ G $ and again $ (\lambda_1, \dots, \lambda_n)_\mathcal{B} \in D $. 

In matrix Lie groups right-invariant vector fields take a simple useful form. In fact, let $G$ be a matrix Lie group. It can be then considered as a Lie subgroup of ${\rm GL}(n,\mathbb{K})$. Moreover, it can be proved that $T_AG$, for any $A\in G$, can be identified with a certain subspace of the space $\mathcal{M}_n(\mathbb{K})$ of $n\times n$ square matrices with coefficients in $\mathbb{K}$.

Since $ R_A : B\in G \mapsto BA\in  G $, then $ R_{A*,e} (M) = M A \in T_A G,$ for all $M \in T_eG$ and $A\in {\rm GL}(n,\mathbb{K})$. As a consequence, if $X^{\rm R}(e)=M$ at the neutral element $e$ of $G$, namely the identity ${\rm I}$ of the matrix Lie group $G$, then
$X^\mathrm{R}(A) = R_{A*,{\rm I}}(X^\mathrm{R}({\rm I})) = R_{A*,{\rm I}}(M) = MA$. It follows that, at any $ A \in G $, every tangent vector $ B \in T_A G $ can be written as $ B = CA $ for a unique $ C \in T_IG$  \cite{MLC,Hall2}. 

Let us describe some basic facts on Lie group actions on manifolds induced by Lie algebras of vector fields. It is known that every finite-dimensional Lie algebra, $V$, of vector fields on a manifold $N$ gives rise to a (local) Lie group action
\begin{equation}\label{Accion}
\varphi : G \times N \to N,
\end{equation}
whose fundamental vector fields are given by the elements of $V$ and $G$ is a connected and simply connected Lie group whose Lie algebra is isomorphic to $V$ (see \cite{Pa_57}). If the vector fields of $V$ are complete, then the Lie group action (\ref{Accion}) is globally defined. Moreover, $G$ does not always need to be simply connected. The Lie group action $\varphi$ will be crucial in the definition of our integrators, since, as can be seen, relates the Lie group $G$ and the manifold $N$, i.e., the manifold where we are going to define the time-evolution of our Lie systems. In fact, Lie group actions like $\varphi$ are employed to reduce the integration of a Lie system on $N$ to obtaining a particular solution of a Lie system on a Lie group \cite{carinena00,lucas11}. Let us show how to obtain $\varphi$ from $V$, which will be of crucial importance in this work.%(\textcolor{blue}{Es necesario que la accion tenga alguna propiedad especial, que sea propia, etc? Alguna otra propiedad a mencionar?}). 

Let us restrict ourselves to an open neighborhood $ U_G $ of the neutral element of $G$, where we can use canonical coordinates of the second-kind related to a basis $ \{ v_1, \dots,v_r \} $ of $\mathfrak{g}$ with opposite structure constants than $X_1,\ldots,X_r$ (see \cite{carinena00,lucas11}). Then, each $ g \in U_G $ can be expressed as
\begin{equation} \label{ccse}
g = \prod_{\alpha = 1}^{r} \exp( \lambda_\alpha v_\alpha ) , 
\end{equation}
for certain uniquely defined parameters $ \lambda_1, \dots, \lambda_r \in \mathbb{R} $. To determine $\varphi$, we first calculate the curves
\begin{equation} \label{curv_accion}
\gamma_x^\alpha : \mathbb{R} \to N : t \mapsto \varphi(\exp(t v_\alpha), x), \qquad \alpha = 1, \dots, r,\qquad \forall x\in N,
\end{equation}
where $\gamma^\alpha_x$ must be the integral curve of $X_\alpha$ starting from the neutral element of $G$ for $\alpha=1,\ldots,r$. Indeed, for any element $ g \in U_G \subset G $ expressed as in \eqref{ccse}, using the intrinsic properties of a Lie group action, 
\begin{equation*}
\varphi(g, x) = \varphi\left( \prod_{\alpha = 1}^{r} \exp( \lambda_\alpha v_\alpha ), x \right)=\varphi(\exp(\lambda_1v_1),\varphi(\exp(\lambda_2 v_2), \ldots,\varphi(\exp(\lambda_r v_r),x)\ldots)),\qquad \forall x\in N,
\end{equation*} 
and the Lie group action is completely defined for every $ g \in U_G \subset G $.

In this work we will deal with some particular matrix Lie groups, starting from
the general linear matrix group $ \mathrm{GL}(n, \mathbb{K}) $, where we stress that $\mathbb{K}$ may be $\mathbb{R}$ or $\mathbb{C}$. % with its bracket defined as $ [A, B] = AB - BA $ for two elements $A,B \in  \mathrm{GL}(n, \mathbb{C})$ \cite[Example 2.10 (a), pg. 31]{lee00}. 
As well known, any closed subgroup of $\mathrm{GL}(n, \mathbb{K})$ is also a matrix Lie group \cite[Theorem 15.29, pg. 392]{lee00}.

\subsection{Generalities on Lie systems} \label{capLie}

Since we are hereafter dealing with non-autonomous systems of differential equations on a manifold $N$ which are locally of the form
$$
\frac{dx^i}{dt}=X^i(t,x),\qquad x\in N,\quad i=1,\ldots,n=\dim N,\qquad t\in \mathbb{R},
$$
we can uniquely relate it to a $t$-dependent vector field, i.e., a $t$-parametric family of standard vector fields on $N$, given by
$$
X(t,x)=\sum_{i=1}^n X^i(t,x)\frac{\partial}{\partial x^i},\qquad x\in N,\quad t\in \mathbb{R},
$$
and conversely. We will hereafter write $\mathfrak{X}_t(N)$ for the space of $t$-dependent vector fields on $N$ \cite{lucas11}. Moreover, $X_t:x\in N\mapsto X(t,x)\in TN$ shall hereafter stand for the vector field induced by $X$ at a particular value of $t\in \mathbb{R}$. The smallest Lie algebra of vector fields (in the sense of inclusion) containing all the vector fields $\{X_t\}_{t\in \mathbb{R}}$ will be denoted by ${\rm Lie}(\{X_t\}_{t\in \mathbb{R}})$. We call ${\rm Lie}(\{X_t\}_{t\in \mathbb{R}})$ the minimal Lie algebra of $X$.

A {Lie system} is a nonautonomous first-order system of ODEs that admits a superposition rule. A {superposition rule} for a system $ X $ on $ N $ (the manifold where $X$ evolves) is a map $ \Phi : N^{m} \times N \to N $ such that the general solution $ x(t) $ of $ X $ can be written as $ x(t) = \Phi(x_{(1)}(t), \dots, x_{(m)}(t); \rho) $, where $ x_{(1)}(t), \dots, $ $ x_{(m)}(t) $ is a generic family of particular solutions and $ \rho $ is a point in $ N $ related to the initial conditions of $X$.

A classic example of Lie system is given by Riccati equations  \cite[Example 3.3]{araujo20}, that is,
\begin{equation} \label{Ric}
\frac{dx}{dt} = b_1(t) + b_2(t) x + b_{12}(t) x^2,\qquad x\in \mathbb{R},\qquad t\in\mathbb{R},
\end{equation} 
with $ b_1(t), b_2(t), b_{12}(t) $ being arbitrary functions of $t$. It is known  that the general solution, $x(t)$, of a Riccati equation can be written as
\begin{equation} \label{SupRiccati}
x(t) = \frac{x_{(2)}(t)(x_{(3)}(t) - x_{(1)}(t)) + \rho x_{(3)}(t) (x_{(1)}(t) - x_{(2)}(t))}{(x_{(3)}(t) - x_{(1)}(t)) + \rho (x_{(1)}(t) - x_{(2)}(t))},
\end{equation}
where $ x_{(1)}(t), x_{(2)}(t), x_{(3)}(t) $ are three different particular solutions of \eqref{Ric} and $ \rho \in \mathbb{R} $ is an arbitrary constant. This implies that every Riccati equation admits a superposition rule $\Phi:\mathbb{R}^3\times \mathbb{R}\rightarrow \mathbb{R}$ of the form
$$
\Phi(x_{(1)},x_{(2)},x_{(3)},\rho)=\frac{x_{(2)}(x_{(3)} - x_{(1)}) + \rho x_{(3)} (x_{(1)} - x_{(2)})}{(x_{(3)} - x_{(1)}) + \rho (x_{(1)} - x_{(2)})}.
$$

The conditions that guarantee the existence of a superposition rule are gathered in the Lie theorem \cite[Theorem 44]{lie93}, which also provides a description of the underlying geometry of a Lie system. This theorem asserts that a first-order system $ X $ on $ N $,
\begin{equation} \label{genLiesys1}
\frac{dx}{dt} = X(t, x), \qquad  x \in N, \qquad  t\in \mathbb{R},\qquad X\in\mathfrak{X}_t(N),
\end{equation}
admits a superposition rule if and only if $ X $ can be written as
\begin{equation} \label{genLiesys}
X(t,x) = \sum_{\alpha = 1}^r b_\alpha(t) X_\alpha(x),\qquad x\in N,\qquad t\in \mathbb{R},
\end{equation}
for a certain family $ b_1(t), \dots, b_r(t) $ of $t$-dependent functions and a family of vector fields $ X_1, \dots, $ $ X_r $ on $N$ that generate an $r$-dimensional Lie algebra of vector fields. Moreover, Lie proved that $\dim V\leq m\dim N$, where $m$ is the number of particular solutions of the associated superposition rule. The relation $\dim V\leq m\dim N$ is called the Lie's condition. The Lie algebra $\langle X_1,\ldots,X_r\rangle$ is called a Vessiot-Guldberg (VG) Lie algebra of the Lie system $X$.

The  $t$-dependent vector field on the real line associated with \eqref{Ric} is $ X_R = b_1(t) Y_0 + b_2(t) Y_1 + b_{12}(t) Y_2 $, where $Y_0,Y_1,Y_2$ are vector fields on $\mathbb{R}$ given by
\begin{equation*}
Y_0 = \frac{\partial}{\partial x}, \qquad Y_1 = x \frac{\partial}{\partial x}, \qquad Y_2 = x^2 \frac{\partial}{\partial x}.
\end{equation*}
Since their commutation relations are
\begin{equation}
    [Y_0,Y_1]=Y_0,\quad [Y_0,Y_2]=2Y_1,\quad [Y_1,Y_2]=Y_2,
\end{equation}
the vector fields $ Y_0, Y_1, Y_2 $ generate a VG Lie algebra $V_R$ isomorphic to $ \mathfrak{sl} (2, \mathbb{R}) $.
Then, Lie theorem guarantees that \eqref{Ric} admits a superposition rule, which is precisely the one shown in \eqref{SupRiccati}. Note that Lie's condition is satisfied since $\dim V_R\leq 3\cdot \dim \mathbb{R}$, where we have used that the superposition rule \eqref{Ric} depends on three particular solutions.

 In practice, the process of determining if $X$ is a Lie system or not relies on taking different vector fields $X_t$ for different values of $t$ and trying to determine the smallest Lie algebra (in the sense of inclusion) containing them by working out all their successive Lie brackets and checking whether all the $X_t$, with $t\in \mathbb{R}$, belong to the obtained family or not. Frequently, one finds two possible outcomes: 

$\bullet$ a) After some calculations,  the successive Lie brackets of some induced vector fields give rise to an infinite family of linearly independent vector fields. This happens in Abel equations of the form (see \cite{lucas11} and references therein)
$$
\frac{dx}{dt}=x^2+a(t)x^3,\qquad x\in \mathbb{R},\qquad t\in \mathbb{R},
$$
for any non-constant $t$-dependent function $a(t)$. In this case, the associated $t$-dependent vector field reads $X_{A}=X_1+a(t)X_2$, where the vector fields $X_1=x^2\partial/\partial x,X_2=x^3\partial/\partial x$ are such that
$$
\left[x^2\frac{\partial}{\partial x},x^k\frac{\partial}{\partial x}\right]=(k-2)x^{k+1}\frac{\partial}{\partial x},\qquad k=3,4,\ldots.
$$
Then, $X_1,X_2$ generate, along with their successive Lie brackets, an infinite family of linearly independent vector fields (over the reals). This shows that the minimal Lie algebra associated with $X_{A}$ is infinite-dimensional and $X_{A}$ is not a Lie system. 

The finite-dimensional Lie algebras of analytic vector fields are also classified on $\mathbb{R}$ and $\mathbb{R}^2$ under quite general conditions \cite{GKO92,Li80,araujo20}. If one associates a differential equation (\ref{genLiesys1}) with its $t$-dependent vector field $X$ and it turns out that the vector fields $\{X_t\}_{t\in \mathbb{R}}$ do not generate a Lie algebra in one of the known types on $N$, then $X$ is not a Lie system. 

$\bullet$ b) In another particular case, one finds that all the $\{X_t\}_{t\in \mathbb{R}}$ are contained in a finite-dimensional Lie algebra of vector fields and $X$ becomes a Lie system.

The Lie theorem yields  that  every Lie system $ X $ is related to (at least) one Vessiot-Guldberg (VG) Lie algebra, $V$, that satisfies that Lie($\{ X_t \}_{t \in \mathbb{R}}$) $\subset V $. This implies that the minimal Lie algebra has to be finite-dimensional \cite{lucas11}. A Lie system may have different VG  Lie algebras. But the minimal Lie algebra is unique.

Let us give a particular illustrative example. Consider the Bernoulli equation
$$
\begin{gathered}
\frac{dx}{dt}=a_1(t)x+a_2(t)x^2,\qquad x\in \mathbb{R},\qquad t\in \mathbb{R},
\end{gathered}
$$
where $a_1(t),a_2(t)$ are arbitrary $t$-dependent functions, 
which is a Lie system as it is related to the $t$-dependent vector field $X_B=a_1(t)Y_1+a_2(t)Y_2$ with $Y_1=x\partial/\partial x$, $Y_2=x^2\partial/\partial x$, and $[Y_1,Y_2]=Y_2$. 

If we consider the particular case $a_1(t)=\cos t$ and $a_2(t)=\sin t$, then the minimal Lie algebra of $X_B$ is the Lie algebra containing all $(X_{B})_t$ with $t\in \mathbb{R}$. In particular, for $t=0$ and $t=\pi/2$, one has $X_0=Y_1$ and $X_{\pi/2}=Y_2$, respectively. Hence, the minimal Lie algebra  of $X$ must contain $Y_1,Y_2$. But $Y_1,Y_2$ span a finite-dimensional Lie algebra and $(X_B)_t$, for every $t\in \mathbb{R}$,  is a linear combination of $Y_1,Y_2$. Hence, $V_B=\langle Y_1,Y_2\rangle$ is the smallest Lie algebra of $X_B$. There exists another VG  Lie algebra of $X_B$ given by
$$
Y_0=\frac{\partial}{\partial x},\qquad Y_1=x\frac{\partial}{\partial x},\qquad Y_2=x^2\frac{\partial}{\partial x}.
$$
The systems related to $\sum_{\alpha=0}^2a_\alpha(t)Y_\alpha$, for arbitrary $t$-dependent functions $a_0(t),a_1(t),a_2(t)$, are the Riccati equations (\ref{Ric}). In other words, using the VG  Lie algebra $V_R=\langle Y_0,Y_1,Y_2\rangle$, one has that  Bernoulli equations with a quadratic term are particular cases of Riccati equations. Riccati equations have then a superposition rule depending on three solutions. Indeed, Lie's condition reads $\dim V_R=1\cdot 3$. Meanwhile, for Bernoulli equations, $\dim V_B=2\cdot 1$, as they are known to have a superposition rule depending on two solutions \cite{lucas11}, but they also accept the superposition rule for Riccati equations, which depends on three particular solutions.

Let us analyze the relevance of using the minimal or a VG Lie algebra of a Lie system. It is known that  $f\in C^\infty(N)$ is a first integral of a Lie system $X$ on $N$ if and only if $f$ is a common first-integral of all the elements of the minimal Lie algebra of the Lie system \cite{araujo20}. Relevantly, it many happen that $f$ is a constant of motion of $X$, but it is not a first integral of all the vector fields of a certain VG  Lie algebra associated with $X$. Then, one may prefer to know the minimal Lie algebra: it characterizes autonomous constants of motion of Lie systems \cite{BCHLS13, araujo20}. Moreover, recall that Lie's condition shows that, if $\dim V$ is larger for a VG Lie algebra $V$ of a Lie system $X$, the superposition rule associated with $X$ depends on more particular solutions, which is not desirable in general. Moreover, many methods to determine superposition rules depend on solving a system of partial differential equations that becomes larger when the VG Lie algebra has larger dimension. That is why, again, one may prefer to determine the minimal Lie algebra for Lie systems. 

On the contrary, the use of larger VG  Lie algebras permits the analysis of properties that are common to more Lie systems. Additionally, some methods to derive superposition rules depend on the nature of VG  Lie algebras. For instance, the coalgebra method \cite{BCHLS13,araujo20} is easier to implement for semi-simple VG Lie algebras, which have the so-called Casimir functions. This implies that it may be simpler to deal with larger VG  Lie algebras if they are semi-simple and have Casimir functions, than with smaller ones that do not have them. 

Let us illustrate the ideas of the previous paragraph. Consider the system of differential equations
$$
\begin{gathered}
\frac{dx_1}{dt}=a(t)x_1+b(t)x_1^2,\\
\frac{dx_2}{dt}=a(t)x_2+b(t)x_2^2,\\
\frac{dx_3}{dt}=a(t)x_3+b(t)x_3^2,\\
\end{gathered}
$$
appearing in the calculus of the superposition rule for Bernoulli equations  with a square term (cf. \cite{lucas11}). In this case, one can relate the above system to the $t$-dependent vector field
$$
Z=a(t)Z_2+b(t)Z_3,\qquad Z_2=\sum_{i=1}^3x_i\frac{\partial}{\partial x_i},\qquad Z_3=\sum_{i=1}^3x^2_i\frac{\partial}{\partial x_i},\qquad [Z_2,Z_3]=Z_3.
$$
Hence, $Z$ is a Lie system. A superposition rule for Bernoulli equations depending on two particular solutions can be obtained by a non-constant first integral for $Z_2,Z_3$ (see \cite{lucas11} for details), which amounts to solving directly a system of PDEs. If one considers Bernoulli equations as particular cases of Riccati equations, the associated superposition rule can be obtained by solving a system of PDEs on $\mathbb{R}^4$ and it depends on three particular solutions \cite{lucas11}. Fortunately, the superposition rule for Riccati equations can be obtained in a simpler manner than for Bernoulli equations via a Casimir function for $\mathfrak{sl}(2,\mathbb{R})$ \cite{araujo20}.

%This minimal Lie algebra is the Vessiot-Guldberg Lie algebra (VG).

%\begin{theorem}[Abbreviated Lie--Scheffers theorem] \label{ALST}
%A system $ X $ is a Lie system if and only if its VG Lie algebra is finite dimensional.
%\end{theorem} 

\subsubsection{Automorphic Lie systems}

It is worth noting that the general solution of a Lie system on $N$ with a VG Lie algebra $V$, can be obtained from a single particular solution of a Lie system on a Lie group $G$ whose Lie algebra is isomorphic to $V$. These are the so-called {automophic Lie systems} \cite[\S 1.4]{lucas11}. As the automorphic Lie system notion is going to be central in our paper, let us study  it in some detail (see \cite{lucas11} for details).

\begin{definition} An {automorphic Lie system} is a $t$-dependent system of first-order differential equations on a Lie group $G$ of the form
\begin{equation}\label{rhshat}
\frac{dg}{dt}=\sum_{\alpha=1}^rb_\alpha(t)X_\alpha^R(g),\qquad g\in G,\quad t\in\mathbb{R},
\end{equation}
where $\{X_1^R,\ldots,X_r^R\}$ is a basis of the space of right-invariant vector fields on $G$ and $b_1(t),\ldots,b_r(t)$ are arbitrary $t$-dependent functions.  Furthermore, we shall refer to the right-hand side of equation \eqref{rhshat} as $\widehat{X}^G_R(t,g)$, i.e., $\widehat{X}^G_R(t,g)=\sum_{\alpha = 1}^{r} b_\alpha(t) X_\alpha^\mathrm{R}(g)$ (see \cite[\S 1.3]{lucas11} for details).
\end{definition}

 Because of right-invariant vector fields, systems in the form of
$ \widehat{X}_R^G $ have the following important property.
\begin{proposition}\label{Prop:RA} 
Given a Lie group $ G $ and a particular solution $ g(t) $ of the Lie system defined on $ G $, as 
\begin{equation}\label{LiesysLiegroup}
\frac{dg}{dt} = \sum_{\alpha = 1}^{r} b_\alpha(t) X_\alpha^\mathrm{R}(g)=\widehat{X}^G_R(t,g), \, 
\end{equation}
where $ b_1(t),\ldots,b_r(t) $ are arbitrary $t$-dependent functions and $ X_1^\mathrm{R},\ldots,X_r^\mathrm{R}$ are right-invariant vector fields, we have that $ g(t) h $ is also a solution of \eqref{LiesysLiegroup} for each $ h \in G $. 
\end{proposition}

An immediate consequence of Proposition \ref{Prop:RA} is that, once we know a particular solution $g(t)$ of  $ \widehat{X}^G_R $, any other solution can be obtained simply by  multiplying $g(t)$ on the right by any element in $ G $. More concretely, the solution $ h(t) $ of (\ref{LiesysLiegroup}) with initial condition $ h(0) = g(0)h_0 $ can be expressed as $h(t)=g(t) h_0 $. This justifies that henceforth we only worry about finding one particular solution $g(t)$ of $ \widehat{X}^G_R $, e.g. the one that fulfills $ g(0) = e $. The previous result can be understood in terms of the Lie theorem or via superposition rules. In fact, since \eqref{LiesysLiegroup} admits a superposition rule $\Phi:(g,h)\in G\times G\mapsto gh\in G$, the system (\ref{Prop:RA})   must be a Lie system. Alternatively, the same result follows from the Lie Theorem and the fact that the right-invariant vector fields on $G$ span a finite-dimensional Lie algebra of vector fields.

There are several reasons to study automorphic Lie systems. One is that they can be locally written around the neutral element of its Lie group in the form 
$$
\frac{dA}{dt}=B(t)A,\qquad A\in {\rm GL}(n,\mathbb{K}),\quad B(t)\in \mathcal{M}_n(\mathbb{K}),
$$ 
where we recall that $\mathcal{M}_n(\mathbb{K})$ is the set of $n\times n$ matrices with coefficients in $\mathbb{K}$.

%We will focus on geometric methods, since when one applies classical numerical methods (linear multi-step or Runge--Kutta) on Lie groups, they generally provide solutions that do not belong in the group, and this is definitely a geometric property that we need to preserve.
The main reason to study automorphic Lie systems is given by the following results, which show how they can be used to solve any Lie system on a manifold. Let us start with a Lie system $ X $ defined on $ N $. Hence, $X$ can be written as
\begin{equation} \label{Lie-Scheffer}
\frac{dx}{dt} = \sum_{\alpha = 1}^{r} b_\alpha(t) X_\alpha(x) ,
\end{equation}
for certain $t$-dependent functions $ b_1(t),\ldots,b_r(t) $ and vector fields $ X_1,\ldots,X_r\in  \mathfrak{X}(N) $ that generate an $r$-dimensional  dimensional VG Lie algebra $V$.
The VG Lie algebra $V$ is always isomorphic to the Lie algebra $ \mathfrak{g} $ of a certain  Lie group $G$. The VG Lie algebra $V$ gives rise to a (local) Lie group action $\varphi:G\times N\rightarrow N$ whose fundamental vector fields are those of $V$. In particular, there exists a basis $\{v_1,\ldots,v_r\}$ of $\mathfrak{g}$ so that
$$
\frac{d}{dt}\bigg|_{t=0}\varphi(\exp(tv_\alpha),x)=X_\alpha(x),\qquad \alpha=1,\ldots,r,\qquad \forall x\in N.
$$
In other words, $\varphi_\alpha:(t,x)\in \mathbb{R}\times N\mapsto \varphi(\exp(tv_\alpha),x)\in N$ is the flow of the vector field $X_\alpha$ for $\alpha=1,\ldots,r$. Note that if $[X_\alpha,X_\beta]=\sum_{\gamma=1}^rc_{\alpha\beta}^\gamma X_\gamma$ for $\alpha,\beta=1,\ldots,r$ and certain constants $c_{\alpha\beta}^\gamma$, then $[v_\alpha,v_\beta]=-\sum_{\gamma=1}^rc_{\alpha\beta}^\gamma v_\gamma$  for $\alpha,\beta=1,\ldots,r$ (cf. \cite{carinena00}).

%Each vector field $X^R_i$ on $ G $ defines a flow of the form $\Phi(t,g)=\exp(t X^R_i(e))$. % As they are defined in terms of right-invariant vector fields, these flows fulfill that
%\[ \Phi^{(i)}(t, g) = \Phi^{(i)}(t, e) g \qquad \forall \, t \in \mathbb{R}, \quad \forall \, g \in G . \] 

%On the other hand, every $ X_i $ on $N$ defines a flow on $N$ that we denote by $ \Phi_i $. 

To determine the exact form of the Lie group action $ \varphi : G \times N \to N $ as in \eqref{curv_accion}, we impose
\begin{equation} \label{def_accion}
\varphi(\exp(\lambda_\alpha v_\alpha), x) = \varphi_\alpha(\lambda_\alpha, x), \qquad  \, \alpha = 1, \dots, r,\qquad \forall x\in N,
\end{equation}
where $ \lambda_1,\ldots,\lambda_r \in \mathbb{R} $. While we stay in a neighborhood $U$ of the origin of $G$, where every element $g\in U$ can be written in the form 
$$
g=\exp(\lambda_1v_1)\cdot\ldots\cdot \exp(\lambda_rv_r),
$$
then the relations \eqref{def_accion} and the properties of $\varphi$ allow us to determine $ \varphi $ on a subset of $G\times N$. If we fix $ x \in N $ and $\alpha$ in (\ref{def_accion}), the right-hand side of the equality turns into the integral curve of the vector field $ X_\alpha $ with the initial condition $x$, this is why (\ref{def_accion}) holds.
% if we differentiate with respect to the parameter (the lambdas) and equal the expression to zero, we get
% \begin{equation} \label{condAcc} 
% \left. \frac{d}{d\lambda_\alpha} \right|_{\lambda_\alpha = 0} \varphi(\exp(\lambda_i v_\alpha), x) = X_\alpha(x) 
% \end{equation} 
% for each $ i = 1, \dots, r $ and fixed $ x \in N $. 

\begin{proposition} \label{prop.accion} (see \cite{carinena00,lucas11} for details) Let $ g(t) $ be a solution to the system 
\begin{equation}\label{LSLG}
\frac{dg}{dt}=\sum_{\alpha=1}^rb_\alpha(t)X_\alpha^R(g),\qquad  t\in \mathbb{R},\quad g\in G,
\end{equation}
where $b_1(t),\ldots,b_r(t)$ are arbitrary $t$-dependent functions.
Then, $ x(t) = \varphi(g(t), x_0) $ is a solution of $ X=\sum_{\alpha=1}^rb_\alpha(t)X_\alpha $ for every $ x_0 \in N $. In particular, if one takes the solution $ g(t) $ that satisfies the initial condition $ g(0) = e $, then $ x(t) $ is the solution of $ X $ such that $ x(0) = x_0 $. 
\end{proposition}

Let us study a particularly relevant form of automorphic Lie systems that will be used hereafter. If $\mathfrak{g}$ is a finite-dimensional Lie algebra, then Ado's theorem \cite{ado47} guarantees that $\mathfrak{g}$ is isomorphic to a matrix Lie algebra $\mathfrak{g}_M$. Let $ \mathcal{V} = \{M_1, \ldots, M_r\} $ be a basis of $ \mathfrak{g}_M\subset\mathcal{M}_n(\mathbb{R})$. %It is worth noting that $M_1,\ldots,M_r$ are chosen so that the Lie algebra structure with the commutator of matrices gives an isomorphism with the Lie algebra of left-invariant vector fields tangent to $\mathfrak{g}_M\subset {\rm Mat}(n)$.
As reviewed in Section \ref{fundaLie},  each  $ M_\alpha $ gives rise to a right-invariant vector field $X^R_\alpha(g)=M_\alpha g$, with $g\in G$, on $ G $. These vector fields have the opposite  commutation relations than the (matrix) elements of the basis. %It is worth noting that the exponential map in ${\rm Mat}(n)$ is given by the standard expression for matrices. %Indeed, for any matrix $ A \in G $ we have
%\[ [X_i^\mathrm{R}, X_j^\mathrm{R}](A) = X_i^\mathrm{R}\left( X_j^\mathrm{R} (A) \right) - X_j^\mathrm{R}\left( X_i^\mathrm{R} (A) \right) = M_i M_j A - M_j M_i A = [M_i, M_j] A , \]
%where we deduce that $ [X_i^\mathrm{R}, X_j^\mathrm{R}] = [M_i, M_j] $ for all $ i, j = 1, \dots, m $. 

In the case of matrix Lie groups, the system  \eqref{LiesysLiegroup} takes a simpler form. Let $ Y(t) $ be the matrix associated with the element $ g(t)\in G $. Using the right invariance property of each $ X_\alpha^\mathrm{R} $, we have that
\begin{equation*}
\frac{dY}{dt} = \sum_{\alpha = 1}^{r} b_\alpha(t) X_\alpha^\mathrm{R}(Y(t)) = \sum_{\alpha = 1}^{r} b_\alpha(t) R_{Y(t)*, e} \left( X_\alpha^\mathrm{R}(e) \right) = \sum_{\alpha = 1}^{r} b_\alpha(t) R_{Y(t)*, e} (M_\alpha) .
\end{equation*}
We can write the last term as
\[ \sum_{\alpha = 1}^{r} b_\alpha(t) R_{Y(t)*, e} (M_\alpha) = \sum_{\alpha = 1}^{r} b_\alpha(t) M_\alpha Y(t) , \]
in such a way that for matrix Lie groups, the system on the Lie group is
\begin{equation} \label{sistemaGM}
\frac{dY}{dt} = A(t) Y(t) , \qquad Y(0) = I , \qquad \text{with} \quad  A(t) = \sum_{\alpha = 1}^{r} b_\alpha(t) M_\alpha ,
\end{equation}
where $ I $ is the identity matrix (which corresponds with the neutral element of the matrix Lie group) and the matrices $ M_\alpha $ form a finite-dimensional Lie algebra, which is anti-isomorphic to the VG Lie algebra of the system (by anti-isomorphic we mean that the systems have the same constants of structure but that they differ in one sign).

There exist various methods to solve system \eqref{LiesysLiegroup} analytically \cite[\S 2.2]{sardon15}, such as the Levi decomposition \cite{levi05} or the theory of reduction of Lie systems \cite[Theorem 2]{carinena01}. In some cases, it is relatively easy to solve it,  as is the case where $ b_1,\ldots,b_r $ are constants. 
Nonetheless, we are interested in a numerical approach, since we will try to solve the automorphic Lie system  with adapted geometric integrators. The solutions on the Lie group can be straightforwardly translated into solutions on the manifold for the Lie system defined on $N$ via  the Lie group action \eqref{Accion}. This is the main idea behind the numerical integrator that we begin to depic in the following 7 step method, which finally will lead us to numerically integrate Lie systems on the manifold $N$, preserving its geometric properties.

%In the case that concerns us, we do not need to map the solutions in the manifold, since we will be working with Lie systems defined on Lie groups. It is possible to obtain the solution on the manifold via the action, but it has been proven that in this case, the power of our numerical approach is not advantageous with respect to other methods.

%Let us introduce our procedure to numerically integrate Lie systems on Lie groups.

%we will enumerate the steps that we will follow along in order to solve a Lie system defined on a Lie group. We call it the 7-step method.

\medskip

\begin{definition}[The 7 step method: Reduction procedure to automorphic Lie systems]
\label{Metodo7Pasos}

The method can be itemized in the following seven steps:
\begin{enumerate}
\item  Given a Lie system on a manifold $N$, we identify a  VG Lie algebra of vector fields admitting a basis $X_1, \dots, X_r$ and associated with the Lie system on $N$.
\item  We look for a matrix Lie algebra  {\rm $\mathfrak{g}$} that is isomorphic to the VG Lie algebra and determine a basis $\{M_1, \dots, M_r \}\subset \mathcal{M}_{n}(\mathbb{R}) $ with the opposite structure constants than the basis $\{X_1,\ldots,X_r\}$.
\item We integrate the vector fields $ X_1,\ldots,X_r $ to obtain their respective flows $ \Phi_\alpha : \mathbb{R} \times N \to N $ with $\alpha=1,\ldots,r$.
\item Using canonical coordinates of the second kind and the previous flows we construct the Lie group action  $ \varphi : G \times N \to N $ using expressions \eqref{def_accion}. 
\item We define an automorphic Lie system $ \widehat{X}^G_R $ on the Lie group $ G $ associated with $ \mathfrak{g} $ as in \eqref{LiesysLiegroup}.
\item We compute the solution $g(t)$ of the system $\widehat{X}^G_R$ that fulfills $ g(0) = e $.  
\item Finally, we recover the general solution for $ X $ on the manifold $N$ by $ x (t)= \varphi(g(t), x_0) $ for $x_0$ being an arbitrary element of $N$. 
\end{enumerate}
\end{definition}

The 7-step method provides a solution of a Lie system on a manifold by means of one particular solution of an automorphic Lie systems on a Lie group and the action \eqref{Accion}. It is important to emphasize that $x(t)$ obtained in the last step of the 7 step method ``lives'' on the manifold $N$, and therefore carries all its geometric properties. In the next section we introduce how the implementation of the 7-step method is carried out numerically. This is accomplished via the Magnus expansion first and the Runge-Kutta-Munthe-Kaas (RKMK) method afterwards.

Let us provide a simple example to illustrate one of the nice features of Lie systems to analyze their properties. 

$\bullet$ Step 1: Consider the nonautonomous system of differential equation on $\mathbb{R}^2$ of the form 
\begin{equation*}%\label{eq:BadRung}
\frac{dx}{dt}=b_1(t)x,\qquad \frac{dy}{dt}=b_2(t)y^k,
\end{equation*}
where $k$ is odd and $b_1(t),b_2(t)$ are arbitrary $t$-dependent functions. The above system is a Lie system since it is associated with the $t$-dependent vector field $X=b_1(t)X_1+b_2(t)X_2$, where
$$
X_1=x\frac{\partial}{\partial x},\qquad X_2=y^k\frac{\partial}{\partial y}
$$
satisfy that $[X_1,X_2]=0$. Then, $X_1,X_2$ span an abelian two-dimensional VG Lie algebra.

$\bullet$ Step 2: A matrix Lie algebra with a basis closing opposite constants of structure than $X_1,X_2$ is given by diagonal $2\times 2$ matrices with the basis
$$
M_1=\left[\begin{array}{cc}1&0\\0&0\end{array}\right],\qquad M_2=\left[\begin{array}{cc}0&0\\0&1\end{array}\right].
$$
 Then, the associated matrix Lie group is spanned by the product of exponentials of $\langle M_1,M_2\rangle$, which gives invertible matrices
$$
G_2=\left\{\left[
\begin{array}{cc}
a&0\\0&b
\end{array}
\right], a,b>0\right\}.
$$
$\bullet$ Step 3: The integration of the vector fields $X_1,X_2$ gives rise to a Lie group action $\varphi:G_2\times \mathbb{R}^2\rightarrow \mathbb{R}^2$.  The integration of the vector fields $X_1,X_2$ for $k\neq 1$ gives the flows
$$
\begin{gathered}
\Phi_1:(t,x,y)\in \mathbb{R}\times \mathbb{R}^2\mapsto (xe^t,y)\in \mathbb{R}^2,\\ \Phi_2:(t,x,y)\in \mathbb{R}\times \mathbb{R}^2\mapsto \left(xe^t,[y^{1-k}+(1-k)t]^{1/(k-1)}\right)\in \mathbb{R}^2,
\end{gathered}
$$
while for $k=1$ the flows read
$$
\Phi_1:(t,x,y)\in \mathbb{R}\times \mathbb{R}^2\mapsto (xe^t,y)\in \mathbb{R}^2,\qquad \Phi_2:(t,x,y)\in \mathbb{R}\times \mathbb{R}^2\mapsto \left(x,ye^t\right)\in \mathbb{R}^2.
$$
It is worth noting that every finite-dimensional Lie algebra of vector fields gives rise to an integrable distribution. 
 It can be added that $X_1,X_2$ span a generalized distribution  $\mathcal{D}$ on $\mathbb{R}^2$ given by
$$
\begin{gathered}
\mathcal{D}_{(x,0)}=\left\langle \frac{\partial}{\partial x}\right\rangle ,\qquad x\neq 0,\qquad \mathcal{D}_{(0,y)}=\left\langle \frac{\partial}{\partial y}\right\rangle,\qquad y\neq 0,\\
 \mathcal{D}_{(x,y)}=T_{(x,y)}\mathbb{R}^2,\qquad xy\neq 0,\qquad x,y\in \mathbb{R}.
\end{gathered}
$$
This generalized distribution is integrable and its strata are given by: a) the point $(0,0)$; the curves b) $x=0,y>0$, c) $x=0,y<0$, d) $y=0,x>0$, e) $y=0,x<0$; and the four connected components of the region of points $(x,y)$ with $xy\neq 0$. The above strata are the orbits of $\varphi$.

$\bullet$ Step 4: Hence, the Lie group action for $k\neq 1$ reads $\varphi:G_2\times \mathbb{R}^2\rightarrow \mathbb{R}^2$ with
\begin{equation}\label{eq:LieGroupAct}
\varphi\left(\left[
\begin{array}{cc}
a&0\\0&b
\end{array}
\right],(x,y)\right)=\left(xa,[y^{1-k}+(1-k)\ln(b)]^{1/(k-1)}\right),\qquad \forall a,b>0,\qquad \forall (x,y)\in \mathbb{R}^2,
\end{equation}
while the canonical coordinates of the second kind read $\lambda_1=\ln a $ and $\lambda_2=\ln b$. For the case $k=1$, we obtain
\begin{equation}\label{eq:LieGroupAct2}
\varphi\left(\left[
\begin{array}{cc}
a&0\\0&b
\end{array}
\right],(x,y)\right)=\left(xa,yb\right),\qquad \forall a,b>0,\qquad \forall (x,y)\in \mathbb{R}^2,
\end{equation}
It is immediate that the fundamental vector fields of the Lie group action $\varphi$ are given by the Lie algebra $V$.

$\bullet$ Step 5: For every possible value of $k$, the automorphic Lie system reads
$$
\frac{d}{dt}\left[\begin{array}{cc}a&0\\0&b\\\end{array}\right]=\left[\begin{array}{cc}b_1(t)&0\\0&b_2(t)\\\end{array}\right]\left[\begin{array}{cc}a&0\\0&b\\\end{array}\right],\qquad a,b>0.
$$
The induced Lie system on the Lie algebra reads
$$
\frac{d}{dt}\left[\begin{array}{cc}\lambda_1&0\\0&\lambda_2\\\end{array}\right]=\left[\begin{array}{cc}b_1(t)&0\\0&b_2(t)\\\end{array}\right]\left[\begin{array}{cc}\lambda_1&0\\0&\lambda_2\\\end{array}\right],\qquad \lambda_1,\lambda_2\in \mathbb{R}.
$$

$\bullet$ Step 6: One can immediately integrate the above system numerically by means of some method to obtain the particular solution starting from the neutral element.

$\bullet$ Step 7: A numerical solution of the automorphic Lie system will give, via the associated Lie group action, an approximate solution of $X$ that will respect the stratification. %Note that any approximate solution of the Lie system on $\mathbb{R}^2$ will have the problem that a very small mistake in determining $y_n$ will cause the solution to diverge strongly from the real behavior of the solution. 

The above example illustrates an interesting fact. The vector fields of the VG  Lie algebra may not give rise to a globally defined Lie group action (see the denominator in (\ref{eq:LieGroupAct})). In any case, the diffeomorphisms on the manifold $\mathbb{R}^2$ induced by the Lie group action, namely the maps $
\varphi_g:(x,y)
\in 
\mathbb{R}^2\rightarrow \varphi(g,(x,y))\in \mathbb{R}^2$ for $g
\in G_2$, will respect the integral manifolds of $\mathcal{D}$. Moreover, the induced automorphic Lie system on $G_2$ can be solved locally around the neutral element to avoid that the action on the elements  of a numerical solution will not be well-defined.

\section{Numerical methods on matrix Lie groups} \label{cap.discre}

This section adapts known numerical methods on Lie groups to study automorphic Lie systems, which are defined by ordinary differential equations defined on Lie groups of the form \eqref{LSLG}.
 For this purpose, we start by reviewing briefly some fundamentals on numerical methods for ordinary differential equations and Lie groups \cite{hairer93, isaacson66, quarteroni07}, and later focus on two specific numerical methods on Lie groups, the Magnus expansion and RKMK methods \cite{iserles99,iserles05,  munthe98, munthe99, Zanna}. We will rely on {one-step methods with fixed time step}. By that we mean that solutions $x(t)$ of a dynamical system
 
\begin{equation}\label{DynSysNumer}
\frac{dx}{dt}=f(t,x),\quad x(a)=x_0,\quad  f\in\mathfrak{X}_t(N),\quad a\in \mathbb{R},
\end{equation}
are approximated by a sequence of points $x_k=x(t_k)\in N$ with $b>a$, $h=(b-a)/\mathcal{N}$, $t_k=a+kh$ for $k=0,\ldots,\mathcal{N}$,  and
\begin{equation}\label{OneStepMethod}
\frac{x_{k+1}-x_k}{h}=f_h(t_k,x_k,x_{k+1}),
\end{equation}
where $\mathcal{N}$ represents the number of steps into which our time interval is divided. We emphasize here that the left hand side of \eqref{OneStepMethod} symbolically represents a proper discretization of a tangent vector on a manifold.  Note that we cannot ``substract'' elements of the manifold and we consider the expression on a local coordinate system. If $N$ is a vector space, the minus sign recovers its usual meaning). We call $h$ the time step, which is fixed, while  $f_h:\mathbb{R}\times N\times N\rightarrow TN$ is a discrete vector field, which is a given approximation of $f$
in \eqref{DynSysNumer}. As usual, we shall denote the local truncation error by $E_h$, where
\begin{equation}\label{LTE}
E_h=||x_{k+1}-x(t_{k+1})|| 
\end{equation}
is the error at the step $k$. Note that for simplicity, $\|\cdot\|$ is a norm on a local coordinated neighborhood \footnote{In matrix Lie groups and Lie algebras, the main cases we are dealing with, problems can be embedded or are defined in a vector space and we can use the restriction to our manifold of a norm associated with the vector space.} in $N$, and we say that the method is of order $r$ if $E_h=\mathcal{O}(h^{r+1})$ for $h\rightarrow 0$, i.e. $\lim_{h\rightarrow0}|E_h/h^{r+1}|<\infty$. Regarding the global error
\[
E_\mathcal{N}=||x_\mathcal{N}-x(b)||,
\]
we shall say that the method is \textit{convergent} of order $r$ if $E_\mathcal{N}=\mathcal{O}(h^{r})$, when $h\rightarrow 0$. As for the simulations, we pick the following norm in order to define the global error, that is
\begin{equation}\label{NormaMaximo}
E_\mathcal{N}=\max_{k=1,\ldots,\mathcal{N}}||x(t_k)-x_k||.
\end{equation}

Our purpose is to numerically solve the initial condition problem for system \eqref{sistemaGM} defined on a matrix Lie group $ G $ of the form
\begin{equation} \label{eq.grupo}
\frac{dY}{dt} = A(t) Y \qquad \text{with} \qquad Y(0) = I, 
\end{equation}
where $ Y \in G $, while $ A(t)$ is a given $t$-dependent matrix taking values in  $\mathfrak{g} \cong T_e G$, and $ I $ is the identity matrix in $G$. That is, we are searching for a discrete sequence $\{Y_k\}_{k=0,\ldots,\mathcal{N}}$ of matrices in $G$. In an open neighborhood of the zero in $T_eG$, the exponential map defines a diffeomorphism onto an open subset $U$ of the neutral element of $G$ and the problem is equivalent on $U$ to searching for a curve
 $ \Omega(t)$ in $\mathfrak{g} $ such that
\begin{equation} \label{eq.omega}
Y(t) = \exp(\Omega(t)) . 
\end{equation}
This ansatz helps us to transform \eqref{eq.grupo}, which is defined in a nonlinear space, into a new problem in a linear space, namely the Lie algebra $\mathfrak{g}\simeq T_eG$. This is expressed in the classical result by Magnus \cite{Magnus}.

\begin{theorem}[Magnus, 1954]

The solution of the matrix Lie group \eqref{eq.grupo} in $G$ can be written for values of $t$ close enough to zero, as $Y(t)=\exp(\Omega(t))$, where $\Omega(t)$ is the solution of the initial value problem
\begin{equation} \label{eq.algebra} 
\frac{d \Omega}{dt} = \operatorname{dexp}^{-1}_{\Omega(t)}(A(t)),\qquad\quad \Omega(0) = {\bf 0} \, ,
\end{equation}
where $ {\bf 0} $ is the zero element in $T_eG$.
\end{theorem}

When we are dealing with matrix Lie groups and Lie algebras, the map $\mbox{dexp}^{-1}$ is given by
\begin{equation} \label{eq.algebraS}
 \mbox{dexp}^{-1}_{\Omega}(H)= \sum_{j = 0}^{\infty} \frac{B_j}{j!} \operatorname{ad}_{\Omega}^j(H), \, 
\end{equation} 
where $\{B_j\}_{j=0,\ldots,\infty}$ are the Bernoulli numbers and $\mbox{ad}_{\Omega}(H)=[\Omega,H]=\Omega\,H-H\,\Omega$ for every $H\in \mathfrak{g}$. The convergence of the series \eqref{eq.algebraS} is ensured as long as a certain convergence condition is satisfied \cite{Magnus}. %$||\Omega||<\pi Magnus no escribe nada de esto$ .

If we try to integrate \eqref{eq.algebra} applying a numerical method directly (note that, now, we could employ one-step methods \eqref{OneStepMethod} safely), $ \Omega(t) $  might sometimes drift too much away from the origin and the exponential map would not work. This would be a problem, since we are assuming that $ \Omega(t) $ stays in a neighborhood of the origin of $\mathfrak{g}$ where the exponential map defines a local diffeomorphism with the Lie group. Since we still do not know how to characterize this neighborhood, it is necessary to adopt  a strategy  that allows us to resolve \eqref{eq.algebra} sufficiently close to the origin. The thing to do is to change the coordinate system in each iteration of the numerical method (Or keepking the time step $h$ small enough, as we shall show when treating the Magnus methods). In the next lines we explain how this is achieved. 

Consider now the restriction of the exponential map given by
\begin{align*} 
\exp : U_{\mathfrak{g}} \subset \mathfrak{g} &\to \exp(U_\mathfrak{g}) \subset G ,\\
A &\mapsto \exp(A) 
\end{align*}
so that this map establishes a diffeomorphism between
an open neighborhood $ U_{\mathfrak{g}} $ around the origin in $ \mathfrak{g} $ and its image. Since the elements of the matrix Lie group are invertible matrices, the map $ U_{\mathfrak{g}} \to \exp(U_\mathfrak{g})Y_0 \subset G : A \mapsto \exp(A) Y_0 $ from $ U_{\mathfrak{g}} \subset \mathfrak{g} $ to the set
\[ \exp(A)Y_0 = \{ Y \in G \ : \exists X\in U_\mathfrak{g}, Y = X Y_0\} \] is also a diffeomorphism. This map gives rise to the so-called {first-order canonical coordinates centered} at $ Y_0 $. 

As well-known, the solutions of \eqref{eq.algebra}  are curves in $\mathfrak{g}$ whose images by the exponential map are solutions to \eqref{eq.grupo}. In particular, the solution $ \Omega^{(0)}(t) $ of system \eqref{eq.grupo} such that $\Omega^{(0)}(0)$ is the zero matrix in $T_IG$, namely ${\bf 0}$,  
corresponds with the solution $Y^{(e)}(t)$ of the system on $G$ such that
$Y^{(e)}(0)=I$. Now, for a certain $ t = t_k $, the solution $ \Omega^{(t_k)}(t)$ in $\mathfrak{g} $ such that $ \Omega^{(t_k)}(t_k) = {\bf 0}$, corresponds with $ Y^{(e)}(t) $
via first-order canonical coordinates centered at
 $ Y^{(e)}(t_k) \in G $, since 
\[ \exp(\Omega^{(t_k)}(t_k)) Y^{(e)}(t_k) = \exp({\bf 0}) Y^{(e)}(t_k) = Y^{(e)}(t_k) , \]
and the existence and uniqueness theorem guarantees $ \exp(\Omega^{(0)}(t)) = \exp(\Omega^{(t_k)}(t)) Y^{(e)}(t_k) $ around $ t_k $. In this way, we can use the curve
 $ \Omega^{(t_k)}(t)$ and the canonical coordinates centered on  $ Y^{(e)}(t_k) $ to obtain values for the solution of \eqref{eq.grupo} in the proximity of $ t = t_k $, instead of using $ \Omega^{(0)}(t) $. Whilst the curve $ \Omega^{(0)}(t) $ could be far from the origin of coordinates for $t_k$,  we know that $ \Omega^{(t_k)}(t) $ will be close, by definition. Applying this idea in each iteration of the numerical method, we are changing the curve in $ \mathfrak{g} $ to obtain the approximate solution of \eqref{eq.grupo} while we stay near the origin (as long as the time step is small enough). 
 
 Thus, what is left is defining proper numerical methods for \eqref{eq.algebra} whose solution, i.e. $\{\Omega_k\}_{k=0,\ldots,\mathcal{N}}$, via the exponential map, provides us with a numerical solution of \eqref{eq.grupo} remaining in $G$. In other words, the general Lie group method defined this way \cite{iserles99, iserles05} can be set by the recursion
 \begin{equation}\label{GeneralLieGroupMethod}
     Y_{k+1}=e^{\Omega_k}\,Y_k,\qquad k=0,1,2,\ldots
 \end{equation}
Next, we introduce two relevant families of numerical methods providing $\{\Omega_k\}_{k=0,\ldots,\mathcal{N}}$.

\subsubsection*{The Magnus method}
Based on the work by Magnus, the Magnus method was introduced in \cite{iserles99,iserles98}. The starting point of this method is to resolve equation \eqref{eq.algebra} by means of the Picard procedure. This method assures that a given sequence of functions converges to the solution of \eqref{eq.algebra} in a small enough neighborhood of $0\in \mathbb{R}$. Operating, one obtains the \textit{Magnus expansion}
\begin{equation} \label{ex.magnus}
\Omega(t) = \sum_{k = 0}^{\infty} H_k(t) , 
\end{equation}
where each $ H_k(t)$ is a linear combination of iterated commutators. The first three terms are given by
\begin{align*}
H_0(t) &= \int_0^t A(\xi_1) d \xi_1 \, , \\
H_1(t) &= -\frac{1}{2} \int_0^t\left[ \int_0^{\xi_1} A(\xi_2) d\xi_2, A(\xi_1) \right] d\xi_1 \, , \\
H_2(t) &= \frac{1}{12} \int_0^t\left[ \int_0^{\xi_1} A(\xi_2) d\xi_2, \left[ \int_0^{\xi_1} A(\xi_2) d\xi_2, A(\xi_1) \right] \right] d\xi_1 \\
&\qquad \qquad + \frac{1}{4} \int_0^t\left[ \int_0^{\xi_1} \left[ \int_0^{\xi_1} A(\xi_2) d\xi_2, A(\xi_1) \right] d\xi_2, A(\xi_1) \right] d\xi_1 \, .
\end{align*}
Note that the Magnus expansion
\eqref{ex.magnus} converges absolutely in a given norm for every $ t \geq 0 $ such that \cite[p. 48]{iserles05}
\[ \int_0^t \| A(\xi) \| d\xi \leq \int_{0}^{2 \pi} \frac{d \xi}{4 + \xi [1 - \cot(\xi / 2)]} \approx 1,086868702. \]

In practice, if we work with the Magnus expansion we need a way to handle the infinite series and calculate the iterated integrals. Iserles and N\o rsett proposed a method based on binary trees \cite{iserles99, iserles98}. In \cite[\S 4.3]{iserles05} we can find a method to truncate the series in such a way that one obtains the desired order of convergence. Similarly, \cite[\S 5]{iserles05} discusses in detail how the iterated integrals can be integrated numerically.  
In our case, for practical reasons we will implement the Magnus method following the guidelines of Blanes, Casas \& Ros \cite{blanes00}, which is based on a Taylor series of $ A(t) $ in \eqref{eq.grupo} around the point $ t = h/2 $ (recall that, in the Lie group and Lie algebra equations we are setting the initial time $t_0=a=0$). With this technique one is able to achieve different orders of convergence. In particular, we will use the second and fourth order convergence methods \cite[\S 3.2]{blanes00}, although one can build up to eighth order methods.

The second-order approximation is
\[ \exp(\Omega(h)) = \exp(h a_0) + \mathcal{O}(h^3) \]
and the forth-order one reads
\[ \exp(\Omega(h)) = \exp\left( h a_0 + \frac{1}{12} h^3 a_2 - \frac{1}{12} h^3 [ a_0, a_1 ] \right)  + \mathcal{O}(h^5) , \]
where $ \Omega(0) = {\bf 0} $ and 
\[ a_i = \frac{1}{i!} \left. \frac{d^i}{dt^i} A(t) \right|_{t = h/2} \qquad i = 0, 1, 2 . \]

As we see from the definition, the first method computes the first and second derivative of matrix $ A(t) $. Applying the coordinate change in each iteration \eqref{GeneralLieGroupMethod}, we can implement it through the following equations:
\begin{align} \label{met.mag2}
Y_{k + 1} = \exp \left[ h A \left( t_k + \frac{h}{2} \right) \right]  Y_k.  \qquad \text{[Order 2]} \\[10pt]
\left. \begin{gathered} \label{met.mag4}
           Y_{k + 1} = \exp\left( h a_0 +  h^3 (a_2 - [ a_0, a_1 ]) \right) Y_k , \\[3pt]
           t_{1/2} = t_k + \frac{h}{2} , \quad a_0 = A(t_{1/2} ) , \quad a_1 = \frac{\dot{A} (t_{1/2})}{12} , \quad a_2 = \frac{\ddot{A} (t_{1/2})}{24},       \end{gathered} \ \right\rbrace \qquad \text{[Order 4]}
\end{align}
where $\dot A(t_0),\ddot A(t_0)$ stand for the first and second derivatives of $A(t)$ in terms of $t$ at $t_0$.  
Note that the convergence order is defined for the Lie group dynamics \eqref{eq.grupo}. That is, when we say that the above methods are convergent of order 2, for instance, that means $E_\mathcal{N}=||Y_\mathcal{N}-Y(b)||=\mathcal{O}(h^2)$, with $h\rightarrow 0$, for a proper Lie matrix norm. It is worth noting that the approximations of $A(t)$, given by derivatives of $A(t)$ and their Lie brackets, belong to the Lie algebra.
Moreover, it is quite apparent in this method that keeping $h$ small enough ensures that $hA\in U_{\mathfrak{g}}$, i.e., the exponential of the Lie algebra element indeed belongs to the Lie group $G$. 
\subsubsection*{The Runge-Kutta-Munthe-Kaas method}

Changing the coordinate system in each step, as explained in previous sections, the classical RK methods applied to Lie groups give rise to the so-called Runge-Kutta-Munthe-Kaas (RKMK) methods \cite{munthe98, munthe99}. The equations that implement the method are
\begin{align*}
\begin{aligned}
\Theta_j  &= h \sum_{l = 1}^{s} a_{jl} F_l , \\
F_j       &= \operatorname{dexp}_{\Theta_j}^{-1}(A(t_k + c_j h)), \\
\Theta    &= h \sum_{l = 1}^{s} b_l F_l , \\
Y_{k + 1} &= \exp(\Theta) Y_k .
\end{aligned} 
\begin{aligned} \left. 
\vphantom{ \begin{aligned}
              \Theta_j  &= h \sum_{l = 1}^{s} a_{jl} F_l , \\
              F_j       &= \operatorname{dexp}_{\Theta_j}^{-1}(A(t_k + c_j h)),
           \end{aligned}
} \quad \right\rbrace \qquad j = 1, \dots, s , \\
\vphantom{ \begin{aligned}
             \Theta    = h \sum_{l = 1}^{s} b_l F_l , \\
             Y_{k + 1} = \exp(\Theta) Y_k .
           \end{aligned}
}
\end{aligned}  
\end{align*}
where the constants $ \{ a_{jl} \}_{j,l = 1}^s $, $ \{ b_l \}_{l = 1}^s $, $ \{ c_j \}_{j = 1}^s $ can be obtained from a Butcher's table \cite[\S 11.8]{quarteroni07} (note that $s$ is the number of stages of the usual RK methods). Apart from this, we have the consistency condition $ \sum_{l = 1}^s b_l = 1 $.
As the equation that we want to solve comes in the shape of an infinite series, it is necessary to study how we evaluate the function $ \operatorname{dexp}_{\Omega(t)}^{-1} $. For this, we need to use truncated series up to a certain order in such a way that the order of convergence of the underlying classical RK is preserved. Moreover, note that the truncated series do always belong to the Lie algebra, as they are given by linear combinations of Lie brackets of elements of the Lie algebra. If the classical RK is of order $p$ and the truncated series of \eqref{eq.algebra} is up to order $j$, such that
$j \geq p - 2 $, then the RKMK method is of order $p$ (see \cite{munthe98,munthe99} and \cite[Theorem 8.5, p. 124]{hairer06}). Again, this convergence order refers to the equation in the Lie group \eqref{eq.grupo}.

Let us now determine the RKMK method associated with the explicit Runge--Kutta whose Butcher's table is
\begin{center}
\begin{tabular}{c|cccc}
$ 0 $   &         &         &         &         \\
$ 1/2 $ & $ 1/2 $ &         &         &         \\
$ 1/2 $ & $ 0 $   & $ 1/2 $ &         &         \\
$ 1 $   & $ 0 $   & $ 0 $   & $ 1 $   &         \\ \hline
        & $ 1/6 $ & $ 1/3 $ & $ 1/3 $ & $ 1/6 $ 
\end{tabular} 
\end{center}
that is a Runge-Kutta of order 4 (RK4). This implies that we need to truncate the series $\operatorname{dexp}_{\Omega(t)}^{-1} $ at $ j = 2 $:
\begin{equation} \label{dexp.2} 
\operatorname{dexp}_{\Omega}^{-1}(A) \approx A - \frac{1}{2} [\Omega, A] + \frac{1}{12} [\Omega, [\Omega, A]] . 
\end{equation}
Then, the RKMK implementation for the given Butcher's table is
\begin{equation} \label{RKMK4}
\left. \begin{aligned}
F_1  &= \operatorname{dexp}_{O_n}^{-1}(A(t_k)) , \\
F_2  &= \operatorname{dexp}_{\frac{1}{2} h F_1}^{-1} \left( A \left( t_k + \frac{1}{2} h \right) \right) , \\
F_3  &= \operatorname{dexp}_{\frac{1}{2} h F_2}^{-1} \left( A \left( t_k + \frac{1}{2} h \right) \right) , \\
F_4  &= \operatorname{dexp}_{h F_3}^{-1}(A(t_k + h))  ,
\end{aligned} \quad \right\rbrace \quad 
\begin{gathered}
\Theta    = \frac{h}{6} (F_1 + 2 F_2 + 2 F_3 + F_4) , \\
Y_{k + 1} = \exp(\Theta) Y_k ,
\end{gathered}
\end{equation}
where $ \operatorname{dexp}^{-1} $ is \eqref{dexp.2}.

It is interesting to note that the method obtained in the previous section using the Magnus expansion \eqref{met.mag2}
can be retrieved by a RKMK method associated with the following Butcher's table
\begin{center}
\begin{tabular}{c|cc}
$ 0 $   &         &       \\
$ 1/2 $ & $ 1/2 $ &       \\ \hline
        & $ 0 $   & $ 1 $ 
\end{tabular} 
\end{center}
Since it is an order 2 method, for the computation of 
$ \operatorname{dexp}^{-1} $ one can use
$ \operatorname{dexp}^{-1}_\Omega (A) \approx A $.

\section{Numerical methods and automorphic Lie systems} \label{metodos}
%As we have already mentioned, by changing the coordinate system in each iteration of the numerical method we will be permanently near the origin of coordinates of the Lie algebra
%\eqref{eq.grupo}, if the step is small. 
%The action over the manifold \ref{pasoNtoG} is defined in terms of the exponential map, therefore it is only well posed for elements near the neutral element of the group. Nevertheless, the solution $ Y(t) $ of \eqref{sistemaGM} can drift away from the neutral element. 

So far, we have established in Procedure \ref{Metodo7Pasos} how to construct an analytical solution of a Lie system $N$ through a particular solution $g(t)$ of an automorphic Lie system on a Lie group $G$ that is based on the integration of the VG Lie algebra associated with the Lie system. More exactly, employing the Lie group action $\varphi$ given in \eqref{Accion}, we can use the particular solution $g(t)$ to obtain the general solution of the Lie system on the manifold $N$. On the other hand, in Section 3 we have reviewed some methods in the literature providing a numerical approximation of the solution of \eqref{eq.grupo} remaining in the Lie group $G$ (which accounts for their most remarkable geometrical property).

Now, let us explain how we combine these two elements to construct our new numerical methods and solve \eqref{Lie-Scheffer}. 
Let $ \varphi $  be the Lie group action constructed using  \eqref{def_accion} and consider the solution of the system \eqref{eq.grupo} such that $ Y(0) = I$. Numerically, we have shown that the solutions of \eqref{eq.grupo} can be provided through the approximations of \eqref{eq.algebraS}, say $\{\Omega_k\}_{k=0,\ldots,\mathcal{N}}$, and \eqref{GeneralLieGroupMethod}, as long as we stay close enough to the origin. As particular examples, we have picked the Magnus and RKMK methods in order to get $\{\Omega_k\}_{k=0,\ldots,\mathcal{N}}$ and, furthermore, the sequence $\{Y_k\}_{k=0,\ldots,\mathcal{N}}$. Next, we establish the scheme providing the numerical solution to automorphic Lie systems.

\begin{definition}\label{MetodosLS}
Let us consider a Lie system 
\begin{equation}\label{LSLG2}
\frac{dx}{dt}=\sum_{\alpha=1}^rb_\alpha(t)X_\alpha(x),\qquad  t\in \mathbb{R},\quad x\in N,
\end{equation}
and let 
$$\frac{dY}{dt}=A(t)Y, \qquad A(t)=\sum_{\alpha=1}^rb_{\alpha}(t)M_{\alpha},
$$
be its associated automorphic Lie system in matrix form close to the neutral element of $G$. We define the numerical solution to the Lie system, i.e., $\{x_k\}_{k=0,\ldots,\mathcal{N}}$, via the algorithm given next.

\begin{algorithm} {\rm Lie systems on Lie groups method}
\label{ForcedAlgorithm}
\begin{algorithmic}[1]
\State {\bf Initial data}: $\mathcal{N}, h,\, A(t),\, Y_0=I ,\, \Omega_0={\bf 0} .$
 \State {\bf Numerically solve } 
 $
 \frac{d\Omega}{dt}= \mbox{dexp}^{-1}_{\Omega}A(t)
 $
 \State{\bf Output} $\{\Omega_k\}_{k=1,\ldots,\mathcal{N}}$
    \For {$k= 1,\ldots, \mathcal{N}-1$} 
    
\[
\begin{split}
Y_{k+1}&=e^{\Omega_k}Y_k,\\
x_{k+1}&=\varphi(Y_{k+1},x_k),
\end{split}
\]
    \EndFor
    \State  {\bf Output:} $(x_1,x_2,...,x_\mathcal{N}).$
%\Statex
\end{algorithmic}
  \end{algorithm}

\end{definition}

At this point, we would like to highlight an interesting geometric feature of this method. On the one hand, the discretization is based on the numerical solution of the automorphic Lie group underlying the Lie system, which, itself, is founded upon the geometric structure of the latter. This numerical solution remains on $G$, i.e., $Y_k\in G$ for all $k$, due to the particular design of the Lie group methods (as long as $h$ is small). Given this, our construction respects as well the geometrical structure of the Lie system, since, in principle, it evolves on a manifold $N$. We observe that the iteration 
\[
x_{k+1}=\varphi(Y_{k+1},x_k)
\]
leads to this preservation, since $x_{k+1}\in N$ as long as $Y_{k+1}\in G$ and $x_k\in N$ (we recall that $\varphi:G\times N\rightarrow N$). Note as well that the direct application of a one-step method \eqref{OneStepMethod} on a general Lie system \eqref{Lie-Scheffer} would destroy this structure, even if applied to an ambient Euclidean space. %(\textcolor{blue}{¿Hay algún sistema de Lie que tenga un $N$ que no sea euclídeo? Esto estaría bien para ponerlo en trabajo futuro: "probar el método en un $N$ no euclideo para demostrar la preservación de la variedad"})

For future reference, in regards of the Lie group methods \eqref{GeneralLieGroupMethod}, we shall refer to \eqref{met.mag2} as Magnus 2, to \eqref{met.mag4} as Magnus 4 and to \eqref{RKMK4} as, simply, RKMK (we recall that the last two methods are order 4 convergent).

\section{Numerical integration on curved spaces}
\label{EjemploDetallado}
This section illustrate the effectiveness of our numerical scheme by applying it to a $(\kappa_1,\kappa_2)$-parametric family of Lie systems on curved spaces (see \cite{HLT17} for details). Naturally, these curved spaces shall play the role of the manifold $N$ where the Lie system evolves. Given that their intrinsic geometry is not trivial, they represent an optimal example of how our  methods are better suited than others for the geometric preservation by our discrete solutions. After introducing a relevant class of Lie systems in curves spaces, Section 5.1 applies the 7 step method and the algorithm in Definition \ref{MetodosLS} to construct the geometry preserving numerical method.

For this, we start by
  considering a two-parametric family of $3D$ real Lie algebras, denoted by $\mathfrak{so}_{\k_1,\k_2}(3)$, which depends on two real parameters,  $\k_1$ and $\k_2$. In the literature these Lie algebras are also known as {CK  Lie algebras}~\cite{  CK2d,Gromovb,Gromova,HOS00,HS02, Ros,Yaglom} or  {quasisimple orthogonal algebras} \cite{casimir}.
  The Lie algebra
$\mathfrak{so}_{\k_1,\k_2}(3)$ admits a basis 
$\{P_1, P_2,J_{12}\}$ with structure constants 
\begin{equation}
  [J_{12},P_1]=P_2,\qquad [J_{12},P_2]=-\k_2 P_1,\qquad [P_1,P_2]=\k_1 J_{12}  . \label{ca}
 \end{equation}

It is possible to rescale the basis of $\mathfrak{so}_{\k_1,\k_2}(3)$ and reducing each parameter $\k_a$ $(a=1,2)$ to either $+1$, 0 or $-1$.
%The vanishment of  any $\k_a$ is   equivalent to applying an {In\"on\"u--Wigner contraction}~\cite{graded}. 
 The Lie algebra $\mathfrak{so}_{\k_1,\k_2}(3)$ admits a representation as a matrix Lie algebra of $3\times 3$ real
matrices $M$ satisfying~\cite{CK2d}
\begin{equation}
M^T \bfI\,+\bfI M=0,\qquad \bfI={\rm diag}(1,\kappa_1,\kappa_1\kappa_2), \qquad \kap=(\k_1,\k_2).
\label{ck2}
\end{equation}
If $\bfI$ is not degenerate, $M$ is called an {indefinite orthogonal Lie algebra} $\mathfrak{so}(p,q)$, where $p$ and $q$ are the number of positive and negative eigenvalues
of $\bfI$. In particular, $\{P_1,P_2,J_{12}\}$ can be identified, respectively, with the matrices
\begin{equation}
P_1=-\k_1 e_{01}+e_{10}, \quad
P_2=-\k_1\k_2 e_{02}+e_{20}, \quad
J_{12}=-\k_2 e_{12}+e_{21} ,
\label{cd}
\end{equation}
where $e_{ij}$ is the $3\times 3$ matrix with a single non-zero  entry 1 at row $i$
and column $j$ with $i,j=0,1,2$. 

The matrix exponential of the elements $\mathfrak{so}_{\kappa_1,\kappa_2}(3)$ generate, by successive matrix multiplications of its elements, the referred to as {CK Lie group} ${\rm SO}_{\kappa_1,\kappa_2}(3)$. In particular, one has the  
  one-parametric subgroups of the CK Lie group ${\rm SO}_{\k_1,\k_2}(3)$ of the form
\begin{equation}\label{ce}
\begin{gathered}
{\rm exp}({\lambda_1 P_1})=\left(\begin{array}{ccc}
\Ck_{\k_1}(\lambda_1)&-\k_1\Sk_{\k_1}(\lambda_1)&0 \cr 
\Sk_{\k_1}(\lambda_1)&\Ck_{\k_1}(\lambda_1)&0\cr 
0&0&1
\end{array}\right) ,\,\,
{\rm exp}({\lambda_2 P_2})=\left(\begin{array}{ccc}
\Ck_{\k_1\k_2}(\lambda_2)&0&-\k_1\k_2\Sk_{\k_1\k_2}(\lambda_2)\cr 
0&1&0\cr 
\Sk_{\k_1\k_2}(\lambda_2)&0&\Ck_{\k_1\k_2}(\lambda_2)
\end{array}\right) ,\\[4pt] 
{\rm exp}({\lambda_3 J_{12}})=\left(\begin{array}{ccc}
1&0&0\cr 
0&\Ck_{\k_2}(\lambda_3)&-\k_2\Sk_{\k_2}(\lambda_3)\cr
0&\Sk_{\k_2}(\lambda_3)&\Ck_{\k_2}(\lambda_3)
\end{array}\right), 
\end{gathered}
\end{equation}
where the so-called {$\k$-dependent cosine} and {sine} functions \cite{CK2d,HOS00, HS02} take the form
\begin{equation}
\Ck_{\k}(\lambda)=\sum_{l=0}^{\infty}(-\k)^l\frac{\lambda^{2l}} 
{(2l)!}=\left\{
\begin{array}{ll}
  \cos {\sqrt{\k}\, \lambda}, &\quad  \k>0, \\ 
\qquad 1,  &\quad
  \k=0, \\
{\rm ch}\, {\sqrt{-\k}\, \lambda}, &\quad   \k<0, 
\end{array}\right.  
\nonumber
\end{equation}
\begin{equation}
   \Sk{_\k}(\lambda) =\sum_{l=0}^{\infty}(-\k)^l\frac{\lambda^{2l+1}}{ (2l+1)!}
= \left\{
\begin{array}{ll}
  \frac{1}{\sqrt{\k}} \sin {\sqrt{\k}\, \lambda}, &\quad  \k>0, \\ 
\qquad \lambda,  &\quad
  \k=0, \\ 
\frac{1}{\sqrt{-\k}} {\rm sh}\, {\sqrt{-\k}\, \lambda}, &\quad  \k<0. 
\end{array}\right.  
\label{ccj}\nonumber
\end{equation}
Similarly to standard trigonometry, the {$\k$-tangent} and the {$\k$-versed sine} (or versine) read
\begin{equation}
\Tk_{\k}(\lambda) =\frac{\Sk_\k(\lambda)}{ \Ck_\k(\lambda)} ,\qquad \Vk_{\k}(\lambda) =\frac 1\k \left(1-\Ck_\k(\lambda) \right)    .
\label{cj}
\end{equation}
These $\k$-functions cover both the usual circular $(\k>0)$ and hyperbolic $(\k<0)$ trigonometric functions. Moreover, $\k$-functions reduce when $\k=0$ to the parabolic functions $\Ck_{0}(\lambda)=1$, $   \Sk{_0}(\lambda) =   \Tk{_0}(\lambda) =\lambda$ and $\Vk_{0}(\lambda) =\lambda^2/2$. 

Some  relations  for the above $\k$-functions read
 \begin{equation}
 \Ck^2_\k(\lambda)+\k\Sk^2_\k(\lambda)=1,\qquad  \Ck_\k(2\lambda)= \Ck^2_\k(\lambda)-\k\Sk^2_\k(\lambda), \qquad \Sk_\k(2\lambda)= 2 \Sk_\k(\lambda) \Ck_\k(\lambda) ,
\label{za}\nonumber
\end{equation}
and their derivatives \cite{HOS00} are given by 

\begin{equation}
 \frac{ {\rm d}}{{\rm d} \lambda}\Ck_\k(\lambda)=-\k\Sk_\k(\lambda),\quad\        \frac{ {\rm d}}
{{\rm d} \lambda}\Sk_\k(\lambda)= \Ck_\k(\lambda)  ,\quad\    
\frac{ {\rm d}}
{{\rm d} \lambda}\Tk_\k(\lambda)=  \frac{1}{\Ck^2_\k(\lambda) } ,\quad\    
\frac{ {\rm d}}
{{\rm d} \lambda}\Vk_\k(\lambda)=  {\Sk_\k(\lambda) }.
\label{zb}
\end{equation}

Define $H_0= {\rm  SO}_{\k_2}(2) $ to be the Lie subgroup of ${\rm SO}_{\kappa_1,\kappa_2}(3)$ obtained by matrix exponentiation of the Lie algebra  ${\mathfrak{h}_0}=\langle J_{12}\rangle$. The CK family of $2D$ homogeneous
spaces  is given by
\begin{equation}
{\mathbf
S}^2_{[\k_1],\k_2} = {\rm  SO}_{\k_1,\k_2}(3)/{\rm  SO}_{\k_2}(2).
\label{cc}
\end{equation}
The (possibly degenerate) metric defined by   ${\bfI}$ in (\ref{ck2}) on $T_e{\rm SO}_{\k_1,\k_2}(3)\simeq \mathfrak{so}_{\k_1,\k_2}(3)$ can be extended to a right-invariant metric on the whole $SO_{\k_1,\k_2}(3)$ by right translation, and then projected onto ${\mathbf
S}^2_{[\k_1],\k_2}$. Then, the CK family becomes a symmetric space relative to the obtained metric. Moreover, $\k_1$ becomes the constant (Gaussian) {\em curvature} of the space, while $\k_2$ determines the {\em signature} of the metric through ${\rm diag}(+,\k_2)$.

The matrix realization induced by (\ref{ce}) leads to the identification of ${\rm SO}_{\k_1,\k_2}(3)$ with the isometries of the bilinear form  $\bfI$. More in detail,
$$
g\in {\rm SO}_{\k_1,\k_2}(3) \Rightarrow g^T \bfI\, g=\bfI,
\label{ck}
$$
which allows one to define a Lie group action of ${\rm SO}_{\k_1,\k_2}(3)$ on $\mathbb{R}^3$ by isometries of $\bfI$.

The subgroup $  {\rm SO}_{\k_2}(2)=\{{\rm exp}({\lambda J_{12}}):\lambda\in \mathbb{R}\} \rangle$ is the isotropy subgroup of the point
$O=(1,0,0)$, which is taken as the {\em origin} in the space
$\mathbf S^2_{[\k_1],\k_2}$. Hence, ${\rm SO}_{\k_1,\k_2}(3)$ becomes an isometry  group of $\mathbf S^2_{[\k_1],\k_2}$.%, so that  $J_{12}$ is a  rotation generator, while $P_1$ and $P_2$ move $O$ along two basic geodesics $l_1$ and $l_2$, which are orthogonal at $O$, so behaving as translation generators.

The orbit of $O$ is contained in the submanifold of $\mathbb{R}^3$ given by
\begin{equation}
\Sigma_\kap =\{v:=(x_0,x_1,x_2)\in \mathbb{R}^3: \bfI(v,v)=\ x_0^2+\k_1  x_1^2+  \k_1\k_2 x_2^2    =1\} .
\label{cl}
\end{equation}
This orbit can be identified with the space  ${\mathbf
S}^2_{[\k_1],\k_2}$. The coordinates $\{x_0,x_1,x_2\}$ on $\mathbb{R}^3$ that
  satisfy  the constraint (\ref{cl}) on $\Sigma_\kap$,  are called {\em ambient}. In these variables, the metric on ${\mathbf
S}^2_{[\k_1],\k_2}$ can be brought into the form
%comes from    the flat ambient metric in $\mathbb R^{3}$ divided by the curvature $\k_1$ and
%restricted to $\Sigma_\kap$, namely
\begin{equation}
{\rm d} s_\kap^2=\left.\frac {1}{\k_1}
\left({\rm d} x_0^2+   \k_1 {\rm d} x_1^2+   \k_1 \k_2{\rm d} x_2^2 
\right)\right|_{\Sigma_\kap}  =    \frac{\k_1\left(x_1{\rm d} x_1 +\k_2 x_2{\rm d} x_2 \right)^2}{1-  \k_1  x_1^2-   \k_1  \k_2 x_2^2 }+  {\rm d} x_1^2+   \k_2{\rm d} x_2^2 .
\label{cm}
\end{equation}
When $\k_1=0$, the manifold $\Sigma_\kap$ has two connected components with $x_0\in \{-1,1\}$, which ensures that ${\rm d}s_\kap^2$ is well-defined.

The ambient coordinates  can be parametrized on $\Sigma_\kap$ through two intrinsic variables in different ways \cite{HLT17,HS02}. Indeed, one may define the so called {\em  geodesic parallel} $\{x,y\}$ or  {\em geodesic polar} $\{r,\phi\}$ coordinates of a point $Q=(x_0,x_1,x_2)$ in $\mathbf S^2_{[\k_1],\k_2}$ that  are defined via the action of the  one-parametric subgroups (\ref{ce}) on $O$~\cite{HS02} given by
$$
(x_0,x_1,x_2)^T = \exp(xP_1) \exp(yP_2)O^T= \exp(\phi J_{12}) \exp(r P_1)O^T ,
\label{cn}
$$
yielding
\begin{equation*}
\begin{gathered}
 x_0=\Ck_{\kk_1}(x)\Ck_{\kk_1\k_2}(y)=\Ck_{\kk_1}(r),\qquad x_1=\Sk_{\kk_1}(x)\Ck_{\kk_1\k_2}(y)=\Sk_{\kk_1}(r)\Ck_{\k_2}(\te) , \\\,\,
 x_2=\Sk_{\kk_1\k_2}(y)=\Sk_{\kk_1}(r)\Sk_{\k_2}(\te) .
 \end{gathered}
%\label{co}
\end{equation*}
Using these relations in the metric (\ref{cm}) and applying (\ref{zb}), one obtains
\begin{equation*}
{\rm d} s_\kap^2=\Ck^2_{\k_1\k_2}(y){\rm d} x^2 + \k_2{\rm d} y^2   =      {\rm d} r^2+\k_2  \Sk^2_{\k_1}(r)  {\rm d} \phi^2 .
%\label{cp}
\end{equation*}

Different values of  $(\k_1,\k_2)$ lead to different spaces. The case $\k_2>0$ gives rise to Riemannian spaces. Any case with $\k_2>0$ and $\k_1<0$ leads to a two-sheeted hyperboloids. Its upper sheet is called $\mathbf {H^2}$, namely the part with $x_0\geq 1$, the {Lobachevsky space}.
Meanwhile, $\k_1= 0$ describes two Euclidean planes $x_0=\pm 1$. We will call the one with $x_0=+1$ Euclidean space $\mathbf {E^2}$.
The instances with $\k_2<0$ define {pseudo-Riemannian spaces or Lorentzian spacetimes}. In this case, a positive Gaussian curvature  $\k_1$ induces a {$(1+1)$D anti-de Sitter spacetime} $\mathbf {AdS^{1+1}}$;  if $\k_1<0$, we find the {$(1+1)$D  de Sitter spacetime $\mathbf {dS^{1+1}}$}; or the flat case with $\k_1=0$, aka the $(1+1)$D  Minkowskian spacetime $\mathbf {M^{1+1}}$. %In all cases for $\k_2<0$, the $J_{12}$, $P_1$, and $P_2$  correspond to the infinitesimal generators of boosts, time translations, and spatial translations, respectively. 
If  $\k_2=0$ $(c=\infty)$, we encounter {Semi-Riemannian spaces or Newtonian spacetimes}, in which the metric  (\ref{cm}) is degenerate and its kernel is an integrable foliation of $\mathbf S^2_{[\k_1],0}$ that is invariant under the 
  action of  the CK group   ${\rm SO}_{\k_1,0}(3)$ on  $\mathbf S^2_{[\k_1],0}$.  There appears a well-defined subsidiary metric  ${\rm d} {s'}^2 ={\rm d} s_\kap^2/\k_2 $  restricted to each leaf, which in the coordinates $(x,y)$  read~\cite{HS02}
  \begin{equation}
  {\rm d} s^2= {\rm d} x^2 ,\qquad   {\rm d}{s'}^2= {\rm d} y^2 \quad {\rm{on}} \quad x= \, {\rm constant} .
  \nonumber
  \end{equation}
For $\k_1>0$ we find the {$(1+1)$D oscillating Newton--Hook (NH) spacetime} $\mathbf {NH_+^{1+1}}$, and for $\k_1<0$ we obtain the {$(1+1)$D expanding NH spacetime} $\mathbf {NH_-^{1+1}}$. The flat space with $\k_1=0$ is just the Galilean $\mathbf {G^{1+1}}$.

%%%%%%%%%%%%%%%%%%%%%%%%%%%%%%%%%%%%%%%%%%%%%%%%%%%%

\subsection{A class of Lie systems on curved spaces}

%We shall hereafter make extensive use of the shorthand notation  $\kap:=(\k_1,\k_2)$. 
Our procedure consists in defining a Lie system $X_{\kap}$ possessing a Vessiot-Guldberg Lie algebra $V_{\kap}$ consisting of infinitesimal symmetries of the metric of the CK space ${\mathbf
S}^2_{[\k_1],\k_2}$. The fundamental vector fields of the Lie group action of ${\rm SO}_{\kap}(3)$ on $\mathbb{R}^3$ by isometries of $\bfI$ are Lie symmetries of ${\rm d}s^2_{\kap}$. Since the action is linear, the fundamental vector fields can be obtained straightforwardly from the $3D$  matrix representation (\ref{cd}).  In ambient coordinates $(x_0,x_1,x_2)$, they read~\cite{HS02},
\begin{equation}
P_1=\k_1 x_1 \frac{\partial}{\partial x_0}   - x_0 \frac{\partial}{\partial x_1} , \qquad
P_2=\k_1\k_2 x_2 \frac{\partial}{\partial x_0}   - x_0 \frac{\partial}{\partial x_2} , \qquad
J_{12}=\k_2 x_2 \frac{\partial}{\partial x_1}   - x_1 \frac{\partial}{\partial x_2}.
\label{da}
\end{equation}

$\bullet$. Step 1. In this case, we are going to use a VG Lie algebra to determine the system we want to analyze. Indeed, the most general Lie system related to $V_{\kap}$ reads%in these ambient coordinates takes the form 
 \begin{equation}\label{tdvf}
  X = b_1(t) P_1 + b_2(t) P_2 + b_{12}(t) J_{12} ,
 \end{equation} 
  where the vector fields $P_1,P_2,J_{12}$ correspond with those in \eqref{da}, and the associated VG Lie algebra has structure constants \eqref{ca}. According to the theory of Lie systems \cite{araujo20}, the integral curves of the time-dependent vector field \eqref{tdvf} are described by the system of ordinary differential equations
\begin{equation}
\left\lbrace
	\begin{aligned}
	\dfrac{dx_0}{dt} &= b_1(t) \kappa_1 x_1 + b_2(t)\kappa_1 \kappa_2  x_2 , \\
	\dfrac{dx_1}{dt} &= -b_1(t) x_0 +  b_{12}(t) \kappa_2x_2 , \\
	\dfrac{dx_2}{dt} &= -b_2(t) x_0 - b_{12}(t) x_1 
	\end{aligned}
\right.
\label{sistema}
\end{equation}
and 
\[ \frac{d(x_0^2 + \kappa_1 x_1^2 + \kappa_1 \kappa_2 x_2^2)}{dt} = 2\left( x_0 \dfrac{dx_0}{dt} + \kappa_1 x_1 \dfrac{dx_1}{dt} + \kappa_1 \kappa_2 x_2 \dfrac{dx_2}{dt} \right)  = 0, \]
which yields that $ I(x_0, x_1, x_2) = x_0^2 + \kappa_1 x_1^2 + \kappa_1 \kappa_2 x_2^2 $ is a constant of motion the Lie system (\ref{sistema}). This invariant will be of utmost importance to show the efficiency of our method when preserving geometric invariants under  numerical integration. Note that determining \eqref{tdvf} for the system \eqref{sistema} can be considered to be the final aim of step 1 of our 7-step method.

$\bullet$ Step 2. It is the moment to consider the following set of matrices
\begin{gather}
M_{P_1} = -\begin{pmatrix}
0 & -\kappa_1 & 0 \\ 1 & 0 & 0 \\ 0 & 0 & 0 
\end{pmatrix}, \qquad
M_{P_2} = -\begin{pmatrix}
0 & 0 & -\kappa_1 \kappa_2 \\ 0 & 0 & 0 \\ 1 & 0 & 0 
\end{pmatrix}, \qquad
M_{J_{12}} = -\begin{pmatrix}
0 & 0 & 0 \\ 0 & 0 & -\kappa_2 \\ 0 & 1 & 0 
\end{pmatrix},
\end{gather}
that have the opposite commutation relations than the vector fields in \eqref{ca}, i.e.,
\[ [M_{P_1}, M_{P_2}] = -\kappa_1 M_{J_{12}}, \qquad [M_{J_{12}}, M_{P_1}] = M_{P_2}, \qquad [M_{J_{12}}, M_{P_2}] = -\kappa_2 M_{P_1} . \]
These are the matrices that we will use as a basis of the Lie algebra of the Lie group to be consider. 

$\bullet$ Step 3. Let us integrate  the vector fields $P_1,P_2,$ and $J_{12}$. To find the flow associated with $ P_1 $, let us solve the system of differential equations $ \{ dx_0/dt = \kappa_1 x_1, \ dx_1/dt = -x_0, \ dx_2/dt = 0 \} $
with initial conditions $ (x_0(0), \ x_1(0),  \ x_2(0))$. The solution is
\begin{equation*}
\left\lbrace \begin{aligned}
x_0(t) &= x_0(0) \operatorname{C_{\kappa_1}}(t) + \kappa_1 x_1(0) \operatorname{S_{\kappa_1}}(t), \\
x_1(t) &= x_1(0) \operatorname{C_{\kappa_1}}(t) -          x_0(0) \operatorname{S_{\kappa_1}}(t), \\
x_2(t) &= x_2(0) ,	
\end{aligned} \right.
\end{equation*}
and the flow $ \Phi_{P_1}:  \mathbb{R} \times \mathbb{R}^3 \to \mathbb{R}^3 $ associated with $ P_1 $ can be expressed in the following manner
\begin{equation}
\Phi_{P_1}(t, (x_0(0), x_1(0), x_2(0))) = (x_0, x_1, x_2), \text{ with } \left\lbrace \begin{aligned}
		x_0 &= x_0(0) \operatorname{C_{\kappa_1}}(t) + \kappa_1 x_1(0) \operatorname{S_{\kappa_1}}(t), \\
		x_1 &= x_1(0) \operatorname{C_{\kappa_1}}(t) -        x_0(0) \operatorname{S_{\kappa_1}}(t), \\
		x_2 &= x_2(0) .	
	\end{aligned} \right.
\end{equation}
Similarly, we calculate the flows for $ P_2 $ and $ J_{12} $ to obtain
\begin{align}
	\Phi_{P_2}(t, (x_0(0), x_1(0), x_2(0))) &= (x_0, x_1, x_2), \text{ with } \left\lbrace \begin{aligned}
		x_0 &= x_0(0) \operatorname{C_{\kappa_1 \kappa_2}}(t) + \kappa_1 \kappa_2 x_2(0) \operatorname{S_{\kappa_1 \kappa_2}}(t), \\
		x_1 &= x_1(0), \\
		x_2 &= x_2(0) \operatorname{C_{\kappa_1 \kappa_2}}(t) -                   x_0(0) \operatorname{S_{\kappa_1 \kappa_2}}(t) ,	
	\end{aligned} \right. \\
	\Phi_{J_{12}}(t, (x_0(0), x_1(0), x_2(0))) &= (x_0, x_1, x_2), \text{ with } \left\lbrace \begin{aligned}
	x_0 &= x_0(0), \\
	x_1 &= x_1(0) \operatorname{C_{\kappa_2}}(t) + \kappa_2 x_2(0) \operatorname{S_{\kappa_2}}(t), \\
	x_2 &= x_2(0) \operatorname{C_{\kappa_2}}(t) -          x_1(0) \operatorname{S_{\kappa_2}}(t) .
\end{aligned} \right.
\end{align}

$\bullet$ Step 4. We now obtain the Lie group action associated with $V_\kap$. Let us briefly review how to obtain the Lie group action related to $V_\kap$. Given an element of the Lie algebra, the canonical coordinates of the second kind permit us to obtain a point in the group (near the origin of coordinates and the neutral element of the algebra, respectively). That is, we have a correspondence between a point in the algebra $ M \in \mathfrak{g} $ determined by the coordinates $ (\lambda_1, \lambda_2, \lambda_3) $, i.e.
\[ M = \lambda_1 M_{P_1} + \lambda_2 M_{P_2} + \lambda_3 M_{J_{12}} \in \mathfrak{g} \]
and the point $ g \in G $ determined by the same coordinates 
\[ g = \exp(\lambda_1 M_{P_1}) \exp(\lambda_2 M_{P_2}) \exp(\lambda_3 M_{J_{12}}) \in G . \]

%Let us determine the Lie group associated with our problem and define canonical coordinates of the second kind on it. 

If we calculate the exponential of $M_{P_1},M_{P_2}, M_{J_{12}}$, we obtain \eqref{ce}.

By multiplying the exponential of these matrices, we obtain the canonical coordinates of the second kind
\begin{multline}
\exp(\lambda_1 M_{P_1}) \exp(\lambda_2 M_{P_2}) \exp(\lambda_3 M_{J_{12}}) = \\
\begin{pmatrix}
\operatorname{C_{\kappa_1}(\lambda_1)} \operatorname{C_{\kappa_1 \kappa_2}(\lambda_2)} & \ast & \ast \\
\operatorname{S_{\kappa_1}(\lambda_1)} \operatorname{C_{\kappa_1 \kappa_2}(\lambda_2)} & \ast & \ast \\
\operatorname{S_{\kappa_1 \kappa_2}(\lambda_2)} & \operatorname{C_{\kappa_1 \kappa_2}(\lambda_2)} \operatorname{S_{\kappa_2}(\lambda_3)} & \operatorname{C_{\kappa_1 \kappa_2}(\lambda_2)} \operatorname{C_{\kappa_2}(\lambda_3)}
\end{pmatrix},
\end{multline}
where we have omitted some matrix entries that are not further needed.  In this way, given a point on the group, we can work out the parameters $ \{ \lambda_1, \lambda_2, \lambda_3 \} $. First, we take the entries $ g_{11} $ y $ g_{21} $ and define $ g = g_{21} / g_{11} $. So,  $ \lambda_1 $ can be expressed as
\begin{equation}
	\lambda_1 = \left\lbrace \begin{array}{ll}
		\dfrac{\arctan g \sqrt{\kappa_1}}{\sqrt{\kappa_1}},                                                       &\text{if } \kappa_1 > 0, \\
		g,                                                                                                        &\text{if } \kappa_1 = 0, \\
		\dfrac{1}{2\sqrt{-\kappa_1}} \log \left( \dfrac{1 + g \sqrt{-\kappa_1}}{1 - g \sqrt{-\kappa_1}} \right),  &\text{if } \kappa_1 < 0 .
	\end{array} \right.
\end{equation}
With the term $ g_{13} $, we can obtain $ \lambda_2 $ as
\begin{equation}
	\lambda_2 = \left\lbrace \begin{array}{ll}
		\dfrac{\arcsin g_{13} \sqrt{\kappa_1 \kappa_2}}{\sqrt{\kappa_1 \kappa_2}},                                                       &\text{if } \kappa_1 \kappa_2 > 0, \\[7pt]
		g_{13}     ,                                                                                                                     &\text{if } \kappa_1 \kappa_2 = 0 ,\\[7pt]
		\dfrac{\log\left( g_{13} \sqrt{-\kappa_1 \kappa_2} + \sqrt{-g_{13}^2 \kappa_1 \kappa_2 + 1} \right)}{\sqrt{-\kappa_1 \kappa_2}}, &\text{if } \kappa_1 \kappa_2 < 0 .
	\end{array} \right.
\end{equation}
And lastly, analogously, defining $ g = g_{32} / g_{33} $ we can obtain $ \lambda_3 $ as
\begin{equation}
\lambda_3 = \left\lbrace \begin{array}{ll}
\dfrac{\arctan g \sqrt{\kappa_2}}{\sqrt{\kappa_2}},                                                       &\text{if } \kappa_2 > 0, \\
g      ,                                                                                                  &\text{if } \kappa_2 = 0, \\
\dfrac{1}{2\sqrt{-\kappa_2}} \log \left( \dfrac{1 + g \sqrt{-\kappa_2}}{1 - g \sqrt{-\kappa_2}} \right),  &\text{if } \kappa_2 < 0 .
\end{array} \right.
\end{equation}

With all of this, by definition, the Lie group action $ \varphi : G \times \mathbb{R}^3 \to \mathbb{R}^3 $ in a point $ g \in G $ and $ \boldsymbol{x}(0) = (x_0(0), x_1(0), x_2(0)) \in \mathbb{R}^3 $ is computed as
\begin{multline*}
\varphi(g, \boldsymbol{x}(0)) = \varphi(\exp(\lambda_1 M_{P_1}) \exp(\lambda_2 M_{P_2}) \exp(\lambda_3 M_{J_{12}}), \boldsymbol{x}(0)) = \\
\varphi(\exp(\lambda_1 M_{P_1}), \varphi(\exp(\lambda_2 M_{P_2}), \varphi( \exp(\lambda_3 M_{J_{12}}), \boldsymbol{x}(0)))) = \\
\Phi_{P_1}(\lambda_1, \Phi_{P_2}(\lambda_2, \Phi_{J_{12}}(\lambda_3, \boldsymbol{x}(0)))).
\end{multline*}

At this point it is interesting to observe that the three vector fields
 $ P_1 $, $ P_2 $ y $ J_{12} $ share the same invariant with the system \eqref{sistema}. This is,   
\begin{equation}
	I(x_0, x_1, x_2) = I\left( \Phi_i (t, (x_0(0), x_1(0), x_2(0)) ) \right) \qquad  \forall t \in \mathbb{R}, \quad \forall (x_0(0), x_1(0), x_2(0)) \in \mathbb{R}^3,
\end{equation}
for each of the flows $\Phi_i$ associated with $\{ P_1, P_2, J_{12} \} $. 
As we have just depicted, the action from the Lie group to the manifold is constructed as the composition of three flows. Moreover, it is apparent that, in spite we are using the ``Euclidean-like'' notation $\mathbb{R}^3$ for our manifold $N$ in this example, it is obvious that it carries a nontrivial geometric structure, as we have already shown. Therefore, our numerical scheme preserves the invariant.

$\bullet$ Steps 5, 6 and 7. Given all these elements, we can implement our numerical scheme. Instead of 
solving \eqref{sistema}, we will solve the following differential equation on the Lie group
\begin{equation}
\frac{dY}{dt} = A(t) Y(t), \qquad Y(0) = I,
\label{sistema_grupo}
\end{equation}
with $ A(t) = b_1(t) M_{P_1} + b_2(t) M_{P_2} + b_{12}(t) M_{J_{12}} $. 

Let us assume that we are in the $ k $-th interaction. This means that we know $ \boldsymbol{x}_k $ and $ Y_k $. To calculate the next point, we apply the numerical scheme \eqref{sistema_grupo} with the initial condition $ Y(0) = Y_k $, obtaining $ Y_{k + 1} $, and being able to compute $ \boldsymbol{x}_{k + 1} = \varphi(Y_{k + 1}, \boldsymbol{x}_k ) $. 

\subsection{Numerical integration of a particular example}
Let us apply our method to \eqref{sistema} with the following coefficients
\[ b_1(t) = t^2, \qquad b_2(t) = \sin t, \qquad b_{12}(t) = \log(t + 1), \]
and constants with values $ (\kappa_1 = 0.8, \kappa_2 = -0.5) $, and initial condition $ \boldsymbol{x}_0 = (1, 1, 1) $, for the interval $ [3, 4] $ and step size $ h = 0.1 $. 

With these parameters our scheme provides the following solution, which is shown overlapped with another solution calculated with a very small step.
We also show the solution obtained with a classical 4th-order Runge-Kutta applied directly to the system, i.e., \eqref{sistema}. 
%One can see this in the graphic depicted the bottom right corner.

\bigskip

\begin{center}
\begin{tikzpicture}
	\begin{axis}[xlabel = {$ t $}, ylabel = {$ x(t) $}, legend pos = outer north east, scale = 0.5, font={\footnotesize}]
		\addplot[line width = 1 pt, color = gray] table[x = t, y = x, col sep = semicolon]{datos/exacta.csv};
		\addplot[mark = o, color = blue]          table[x = t, y = x, col sep = semicolon]{datos/RKMK4.csv};
		\addplot[mark = x, color = red]           table[x = t, y = x, col sep = semicolon]{datos/RK4.csv};
	\end{axis}
\end{tikzpicture}
\begin{tikzpicture}
	\begin{axis}[xlabel = {$ t $}, ylabel = {$ y(t) $}, legend pos = outer north east, scale = 0.5, font={\footnotesize}]
		\addplot[line width = 1 pt, color = gray] table[x = t, y = y, col sep = semicolon]{datos/exacta.csv};
		\addplot[mark = o, color = blue]          table[x = t, y = y, col sep = semicolon]{datos/RKMK4.csv};
		\addplot[mark = x, color = red]           table[x = t, y = y, col sep = semicolon]{datos/RK4.csv};
		\legend{exact, RKMK, RK4};
	\end{axis}
\end{tikzpicture}

\begin{tikzpicture}
	\begin{axis}[xlabel = {$ t $}, ylabel = {$ z(t) $}, legend pos = outer north east, scale = 0.5, font={\footnotesize}]
		\addplot[line width = 1 pt, color = gray] table[x = t, y = z, col sep = semicolon]{datos/exacta.csv};
		\addplot[mark = o, color = blue]          table[x = t, y = z, col sep = semicolon]{datos/RKMK4.csv};
		\addplot[mark = x, color = red]           table[x = t, y = z, col sep = semicolon]{datos/RK4.csv};
	\end{axis}
\end{tikzpicture}
\begin{tikzpicture}
	\begin{axis}[xlabel = {$ t $}, ylabel = {$ I(x, y, z) $}, legend pos = outer north east, scale = 0.5, font={\footnotesize}]
		\addplot[line width = 1 pt, color = gray]    table[x = t, y = exacta, col sep = semicolon]{datos/invariante.csv};
		\addplot[mark = o, color = blue, only marks] table[x = t, y = RKMK,   col sep = semicolon]{datos/invariante.csv};
		\addplot[mark = x, color = red]              table[x = t, y = RK4,    col sep = semicolon]{datos/invariante.csv};
		\legend{exact, RKMK, RK4};
	\end{axis}
\end{tikzpicture}
\end{center}

As we observe in the bottom right plot, the geometric quantity  $ I(x, y, z) $ is exactly preserved, which is not the case with a classical 4th-order RK scheme applied directly to the original Lie system. Once again, this geometric preservation is achieved by means of the specific design of the integrator, which is the main point of our work. Finally, we show a convergence plot (in logarithmic scale) of our scheme for the $x$ component. We employ the definition of the global error $E_{\mathcal{N}}$ given in \eqref{NormaMaximo}, where it is enough to consider an Euclidean norm, since we are in $\mathbb{R}^3$ (with nontrivial curvature). We observe convergence, as in the case of the other two components.

\begin{figure}[ht]
  \centering
  \includegraphics[width=0.4\textwidth]{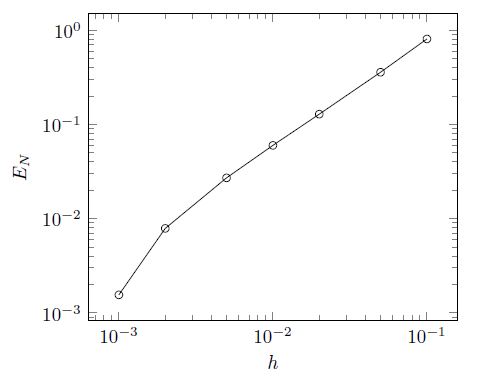} 
\end{figure}

\section{Initial conditions, long term behavior, and numerical methods}\label{Sec:IniLongNum}
To observe more critically when and how our methods are more appropriate than other ones, let us consider  %we will recall the system in \eqref{eq:BadRung}. Note that the time dependent vector field is such that $X_1$ is tangent to the line $x=0$ while $X_2$ is so for the line $y=0$. An initial condition in any of these lines will make the particular solution to remain inside of the the four regions in the plane. Therefore, the problem is very sensitive to initial conditions and numerical methods will produce bad results in this case. But ours will not. This is due to the fact that we respect the integral submanifolds of the generalized distribution spanned by $X_1,X_2$. To illustrate our point, consider the numerical approximation for $b_1(t)=-(1+t^2)$, $b_2(t)=e^t$, $k=1$ and the particular solution with initial condition $(1,0)$. Then, the solution reads as follows
 the nonautonomous system of differential equation on $\mathbb{R}^2$ of the form 
\begin{equation}\label{eq:BadRung2}
\frac{dx}{dt}=b_1(t)y+b_2(t)(x^2+y^2-1)x,\qquad \frac{dy}{dt}=-b_1(t)x+b_2(t)(x^2+y^2-1)y.
\end{equation}
where $b_1(t),b_2(t)$ are arbitrary $t$-dependent functions. What we are going to show is that our methods, which are adapted to the geometric features of (\ref{eq:BadRung2}), more accurately reflect  its long-term behavior and the dependence on initial conditions of particular solutions.

System (\ref{eq:BadRung2}) is a Lie system since it is associated with the $t$-dependent vector field $X=b_1(t)X_1+b_2(t)X_2$, where
$$
X_1=y\frac{\partial}{\partial x}-x\frac{\partial}{\partial y},\qquad X_2=(x^2+y^2-1)\left(x\frac{\partial}{\partial x}+y\frac{\partial}{\partial y}\right)
$$
satisfy that $[X_1,X_2]=0$. Then, $X_1,X_2$ span an abelian two-dimensional VG Lie algebra. Additionally, $X_1,X_2$ span a generalized distribution  $\mathcal{D}$ on $\mathbb{R}^2$ given by
$$
\begin{gathered}
\mathcal{D}_{(x,y):x^2+y^2=1}=\left\langle y\frac{\partial}{\partial x}-x\frac{\partial}{\partial y}\right\rangle ,\qquad \mathcal{D}_{(x,y):x^2+y^2\neq \{1,0\}}=T_{(x,y)}\mathbb{R}^2,\qquad \mathcal{D}_{(0,0)}=0.\\
\end{gathered}
$$
Recall that every finite-dimensional Lie algebra of vector fields gives rise to an involutive distribution.
 This generalized distribution is also integrable and its strata are given by the circle $x^2+y^2=1$, the point $(0,0)$, and the regions $0<x^2+y^2<1$ and $x^2+y^2>1$. 

The integration of the vector fields $X_1,X_2$ gives rise to a Lie group action that respects the leaves of $\mathcal{D}$. The matrix Lie algebra with a basis closing opposite constants of structure than $X_1,X_2$ is spanned by the basis
$$
M_1=\left[\begin{array}{cc}1&0\\0&0\end{array}\right],\qquad M_2=\left[\begin{array}{cc}0&0\\0&1\end{array}\right].
$$
This is the matrix Lie algebra of diagonal $2\times 2$ matrices. Then, its matrix Lie group is spanned by the product of an arbitrary number of  exponentials of elements in $\langle M_1,M_2\rangle$, which gives the Lie group of invertible matrices
$$
G_2=\left\{\left[
\begin{array}{cc}
a&0\\0&b
\end{array}
\right]: a,b>0\right\}.
$$
The integration of the vector fields $X_1,X_2$ leads to the flows
$$
\begin{gathered}
\Phi_1:(t,x,y)\in \mathbb{R}\times \mathbb{R}^2\mapsto (x\cos t+y\sin t,-x\sin t+y\cos t)\in \mathbb{R}^2,\\ \Phi_2:(t,x,y)\in \mathbb{R}\times \mathbb{R}^2\mapsto \frac{(x,y)}{\sqrt{(x^2+y^2-(x^2+y^2-1)e^{2t}}}\in \mathbb{R}^2.
\end{gathered}
$$
And the associated Lie group action reads $\varphi:G_2\times \mathbb{R}^2\rightarrow \mathbb{R}^2$ is such that $(0,0)$ is left invariant by $G_2$, while
\begin{equation}\label{eq:LieGroupAct3}
\varphi\left(\left[
\begin{array}{cc}
a&0\\0&b
\end{array}
\right],(x,y)\right)=\frac{(x\cos \ln a+y\sin \ln a,-x\sin \ln a+y\cos \ln a)}{\sqrt{x^2+y^2-(x^2+y^2-1)b^2)}},
\end{equation}
where the canonical coordinates of the second kind read $\lambda_1=\ln a $ and $\lambda_2=\ln b$. 
This action leaves invariant the circle $x^2+y^2=1$ and the point $x=y=0$. It is important to remark that the Lie group action is local and, given $(x,y)$, the action is only defined for the values of $a,b$ that do not make the denominator in (\ref{eq:LieGroupAct3}) to be negative or zero. In fact, when $b$ is such that the denominator tends to zero, the image of the Lie group action tends to infinite. 

Moreover, the form of $X_1,X_2$ and the associated distribution $\mathcal{D}$ shows that the solutions of (\ref{eq:BadRung2}) must remain in one the strata of $\mathcal{D}$. This will be respected by our methods, but the Runge-Kutta and other numeric methods do need to do so: the circular solution does not need to be respected and solutions may change from one orbit of $\varphi$ to another, which violates one of the fundamental features of (\ref{eq:BadRung2}). This clearly shows that generic numerical methods will not respect geometric features of differential equations, which can be extremely important, for instance, when the particular solutions depend strongly on the initial conditions.

The automorphic Lie system associated with \eqref{eq:BadRung2} reads
$$
\frac{d}{dt}\left[\begin{array}{cc}a&0\\0&b\\\end{array}\right]=\left[\begin{array}{cc}b_1(t)&0\\0&b_2(t)\\\end{array}\right]\left[\begin{array}{cc}a&0\\0&b\\\end{array}\right],\qquad a,b>0.
$$
and the induced Lie system on the abelian Lie algebra $\mathbb{R}^2$ becomes
$$
\frac{d}{dt}\left[\begin{array}{cc}\lambda_1&0\\0&\lambda_2\\\end{array}\right]=\left[\begin{array}{cc}b_1(t)&0\\0&b_2(t)\\\end{array}\right]\left[\begin{array}{cc}\lambda_1&0\\0&\lambda_2\\\end{array}\right],\qquad \lambda_1,\lambda_2\in \mathbb{R}.
$$
Then, the numerical solution of the automorphic Lie system will give, via the associated Lie group action $\varphi$, an approximate solution of $X$ that will respect the stratification.

To verify the usefulness of our method, let us analyse the initial problem with the initial conditions $(0,1)$ and the $t$-dependent coefficients $b_1(t)=(1+t^2)$ and $b_2(t)=e^t$. Note that a Runge--Kutta method will show that a solution with initial condition in $x^2+y^2=1$ will finally escape from that circle and make it to move far away from the real solution.

\begin{center}
\begin{tikzpicture}
	\begin{axis}[xlabel = {$ x $}, ylabel = {$ y $}, legend pos = outer north east]
		\addplot[line width = 1.5 pt, color = gray] table[x = xreal,    y = yreal,    col sep = semicolon]{datos/bifurcacion.csv};
		\addplot[mark = o, color = blue]            table[x = X0.02,    y = Y0.02,    col sep = semicolon]{datos/bifurcacion.csv};
		\addplot[mark = x, color = red]             table[x = X0.01,    y = Y0.01,    col sep = semicolon]{datos/bifurcacion.csv};
        \addplot[mark = *, color = green]           table[x = xRKMK0.1, y = yRKMK0.1, col sep = semicolon]{datos/bifurcacion.csv};
  		\legend{exact, RK ($ h = 0.02 $), RK ($ h = 0.01 $), RKMK ($ h = 0.1 $)};
	\end{axis}
\end{tikzpicture}
\end{center}

\section{Conclusions and outlook}\label{Sec:ConOut}

Given the wide range of spaces, and geometries, where the mathematical and physical dynamical systems evolve, it is always worth to take care of its intrinsic properties when passing to the ``discrete'' side in order to obtain an approximate solution. As extensively showed in the literature, this results on some computational and dynamical benefits. This is the spirit of geometric integration, and the one we uphold in this article, where we take advantage of the geometric structure of Lie systems in order to propose a 7-step method to analytically solve them (mainly, the possibility to reduce such systems to equivalent ones on a Lie group), plus a geometric numeric integrator. We have proven its geometric properties with a wide class of Lie systems evolving on curved spaces. We have also studied examples whose particular solutions may strongly depend on the initial conditions and our methods are more appropriate to study the short- and long-term behavior than other numerical non-geometrical methods.

As for future work, it shall be worth wondering about the numerical features of the integrator, such as consistency and convergence, besides finding new examples which may be of interest in mathematics, physics or other applied sciences. Moreover, it is interesting that linear systems of differential equations have been used to study differential equations close to equilibrium points. Such systems are Lie systems. It will be interesting to use approximations up to second order of differential equations close to equilibrium points, which may potentially lead to Lie systems that reflect more appropriately the properties of the systems under study and the use of the methods of this work (cf. \cite{araujo20} for Lie systems on the plane). There is no classification of Lie algebras of smooth vector fields on the plane whose elements vanish at a certain point (cf. \cite{GKO92}). This is a very interesting topic from the point of view of stability analysis and numerical methods for Lie systems to be developed in a further work.

\section*{Acknowledgements}

J. de Lucas acknowledges a Simons--CRM professorship funded by the Simons Foundation and the Centre de Recherches Math\'ematiques (CRM) of the Universit\'e de Montr\'eal. J. de Lucas would like to thank for the hospitality shown by the members and staff of the CRM during his stay. We also would like to thank the anonymous referees for their suggestions, which allowed us to significantly improve the quality, relevance, and clarity of this work.

 \section*{Data availability}
 The datasets generated during and/or analysed during the current study are available from the corresponding author on reasonable request.


\begin{thebibliography}{60}







\bibitem{Mahony1}
P.A. Absil, R. Mahony, and R. Sepulchre,  {\it Riemannian Geometry of Grassmann Manifolds with a View on Algorithmic Computation}, Acta Appl. Math. {\bf 80}:199-220, 2004.


\bibitem{ado47} 
I.D. Ado, \textit{The representation of Lie algebras by matrices}, Uspekhi Matematicheskikh Nauk \textbf{2}:159-173, 1947.

\bibitem{angelo05} 
R.M. Angelo and W.F. Wreszi\'nski, \textit{Two-level quantum dynamics, integrability and unitary NOT gates}, Phys. Rev. A \textbf{72}:034105, 2005.

\bibitem{BCHLS13}
A. Ballesteros, J.F. Cari\~nena, F.J. Herranz, J. de Lucas, C. Sardón, {\it From constants of motion to superposition rules for Lie-Hamilton systems}, J. Phys. A {\bf 46}:285203, 2013.

\bibitem{blanes06}
S. Blanes, F. Casas, J.A. Oteo and J. Ros, {\it  The Magnus expansion and some of its application}, Phys. Rep. {\bf 470}:5-6, 2006.

\bibitem{blanes00} 
S. Blanes, F. Casas and J. Ros, \textit{Improved high order integrators based on the Magnus expansion}, BIT Nume. Math. \textbf{40}:434-450, 2000.


\bibitem{Moan}
S. Blanes and P.C. Moan, {\it Practical Symplectic Partitioned Runge-Kutta and Runge-Kutta-Nyström Methods}, J. Computational  Appl. Math. {\bf 142}:313-330, 2002.

\bibitem{blasco15} A. Blasco, F.J. Herranz, J. de Lucas and C. Sard\'on, \textit{Lie-Hamilton systems on the plane: Applications and superposition rules}, J. Phys. A  \textbf{48}:345202, 2015.

\bibitem{carinena00} J.F. Cariñena, J. Grabowski and G. Marmo, {\sl Lie-Scheffers Systems: a Geometric Approach}, Napoli Series in Physics and Astrophysics, Bibliopolis, Naples, 2000.

\bibitem{carinena07} J.F. Cariñena, J. Grabowski and G. Marmo, \textit{Superposition rules, Lie theorem and partial differential equations}, Rep. Math. Phys. \textbf{60}:237-258, 2007.

\bibitem{carinena01} J.F. Cariñena, J. Grabowski and A. Ramos, \textit{Reduction of $t$-dependent systems admitting a superposition principle}, Acta Appl. Math. \textbf{66}:67-87, 2001.

\bibitem{carinena09} J.F. Cariñena and J. de Lucas, \textit{Applications of Lie systems in dissipative Milne-Pinney equations}, Int. J. Geom. Meth.  Modern Phys. \textbf{6}:683-699, 2009.

\bibitem{lucas11} 
J.F. Cariñena and J. de Lucas, {\sl Lie Systems: Theory, Generalisations, and Applications}, Diss. Mathematicae \textbf{479}, Warsaw, 2011.

\bibitem{carinena11} 
J.F. Cariñena, J. de Lucas and C. Sardón, \textit{A new Lie systems approach to second-order Riccati equations}, Int. J.  Geom. Meth.  Modern Phys. \textbf{9}:1260007, 2011.

\bibitem{carinena99} J.F. Cariñena and A. Ramos, \textit{Integrability of the Riccati equation from a group theoretical viewpoint}, Int. J. Mod. Phys. A \textbf{14}:1935-1951, 1999.


\bibitem{CortesMartinez} J. Cort\'es and S. Mart\'inez, \textit{Non-holonomic integrators}, Nonlinearity \textbf{14}:1365-1392, 2001.

\bibitem{MLC}
 M.L. Curtis, {\sl Matrix groups}, Springer,  New York, 1984.
 
\bibitem{dominguez06} S. Domínguez, P. Campoy, J.M. Sebastián and A. Jiménez,  {\sl Control en el Espacio de Estado}, Pearson, Educación, 2006.



\bibitem{CK2d}  
A. Ballesteros, F.J. Herranz, M.A. del Olmo, and M.  Santander,  
{\it Quantum structure of the motion groups of the two-dimensional Cayley--Klein geometries},
J. Phys. A {\bf 26}:5801-5823, 1993.


\bibitem{GKO92}
A. Gonz\'alez, N. Kamran, and P.J. Olver, {\it Lie algebras of vector fields on the plane}, Proc. London Math. Soc. {\bf 64}:339-368, 1992.


\bibitem{Gromovb}   
N.A. Gromov, {\em Contractions and analytical continuations of the
classical groups. Unified approach},  Komi  Scienfic Center, Syktyvkar, 1992 (in Russian).

\bibitem{Gromova}
N.A. Gromov and V.I. Man'ko,
{\it  The Jordan--Schwinger representations of
Cayley--Klein groups. I. The orthogonal groups}, 
J. Math. Phys. {\bf 31}:1047-1053, 1990.


\bibitem{hairer06} E. Hairer, C. Lubich and G. Wanner, {\sl Geometric Numerical Integration}, Springer-Verlag, Berlin-Heidelberg, 2006.

\bibitem{hairer93} E. Hairer, S.P. N\o rsett and G. Wanner, {\sl Solving Ordinary Differential Equations I: Nonstiff Problems}, Springer-Verlag, Berlin, 1993.

%\bibitem{Hall}
% B. Hall, {\it Matrix Lie Groups} in {\sl Lie Groups, Lie Algebras, and Representations: An Elementary Introduction}, Springer, Cham, 2015, pp. 3-30.

\bibitem{Hall2}
B.C. Hall, {\sl Lie Groups, Lie Algebras, and Representations: An Elementary Introduction}, Springer, Cham, 2015.

\bibitem{HWA83}
J. Harnad, P. Winternitz and R.L. Anderson, 
{\it Superposition principles for matrix Riccati equations}, 
J. Math. Phys. {\bf 24}:1062, 1983.

\bibitem{hartshorne67} 
R. Hartshorne, {\sl Foundations of Projective Geometry}, W.A. Benjamin, Inc., New York, 1967.


\bibitem{HLT17}
 F.J. Herranz, J. de Lucas and M. Tobolski, {\it 
Lie-Hamilton systems on curved spaces: A geometrical approach},
J. Phys. A {\bf 50}:495201, 2017. 




 
\bibitem{graded}
F.J. Herranz, M. de Montigny, M.A. del Olmo  and M. Santander,
{\it Cayley--Klein algebras as graded contractions of $so(N+1)$}, {J. Phys. A} {\bf 27}:2515--2526, 1994.


\bibitem{HOS00}
F.J. Herranz, R. Ortega and M. Santander,
{\it Trigonometry of spacetimes: a new self-dual approach to a curvature/signature (in)dependent trigonometry},
 J. Phys. A {\bf 33}:4525-4551, 2000.
 
\bibitem{casimir}
F.J. Herranz and M. Santander,
{\it Casimir invariants for the complete family of quasisimple orthogonal algebras},
 {J. Phys. A} {\bf 30}:5411-5426, 1997.
\bibitem{HS02}
F.J. Herranz and M. Santander, 
{\it Conformal symmetries of spacetimes},
J. Phys. A {\bf 35}:6601-6618, 2002.


\bibitem{hussin90} V. Hussin, J. Beckers, L. Gagnon and P. Winternitz, \textit{Superposition formulas for nonlinear superequations}, J. Math. Phys. \textbf{31}:2528-2534, 1990.

\bibitem{iserles99} 
A. Iserles and S.P. N\o rsett, \textit{On the solution of linear differential equations in Lie groups}, Phil. Trans Royal Soc. A \textbf{357}:983-1020, 1999.

\bibitem{iserles98} 
A. Iserles, S.P. N\o rsett and A.F. Rasmussen, \textit{$t$-symmetry and high-order Magnus methods}, Technical Report 1998/NA06, DAMTP, University of Cambridge, 1998.

\bibitem{isaacson66} E. Isaacson and H.B. Keller, {\sl Analysis of Numerical Methods}, John Wiley \& Sons, New York-London-Sydney, 1966.

\bibitem{iserles05} A. Iserles, H. Munthe-Kaas, S. N\o rsett and A. Zanna, {\it Lie-group methods}, Acta Num. {\bf 9}:215-365, 2000.

\bibitem{Ku73}
V. Ku\v{c}era, 
{\it A Review of the Matrix Riccati Equation},
Kybernetika {\bf 9}:42-61, 1973.

\bibitem{Ko83}
L. K\"onigsberger, {\it \"Uber die einer beliebigen differentialgleichung erster Ordnung angeh\"origen selb-st\"andigen Transcendenten}, Acta Math.{\bf 3}:1-48, 1883. 

\bibitem{LL19}
 J. Lange and J. de Lucas, {\it Geometric models for Lie--Hamilton systems on $\mathbb{R}^2$},
 Mathematics {\bf 2019}:7, 1053.
 
\bibitem{lazaro09} 
J.A. Lázaro-Camí and J.P. Ortega, \textit{Superposition rules and stochastic Lie-Scheffers systems}, Ann. Inst. H. Poincaré Probab. Stat. \textbf{45}:910-931, 2009.

\bibitem{lee00} J.M. Lee, {\sl Introduction to Smooth Manifolds}, Graduate Texts in Mathematics {\bf 218}, Springer-Verlag, New York, 2003. 

\bibitem{levi05} 
E.E. Levi, \textit{Sulla struttura dei gruppi finiti e continui}, Atti della Reale Accademia delle Scienze di Torino, 1905.

\bibitem{Li80}
S. Lie, {\sl Sophus Lie's 1880 Transformation group paper},  (Translated by M. Ackerman and comments by R. Hermann), Lie Groups: History, frontiers and applications I, Math. Sci. Press, Mass, 1975. 

\bibitem{lie93} S. Lie and G. Scheffers, {\sl Vorlesungen über continuierliche Gruppen mit geometrischen und anderen Anwendungen}, Teubner, Leipzig, 1893.

\bibitem{LG17}
J. de Lucas and A.M. Grundland,
{\it A Lie systems approach to the Riccati hierarchy and partial differential equations}, J.  Differential Equations {\bf 263}:299-337, 2017.


 
\bibitem{araujo20} 
J. de Lucas and C. Sardón, {\sl A Guide to Lie Systems with Compatible Geometric Structures}, World Scientific, Singapore, 2020. 

\bibitem{Magnus}
W. Magnus, \textit{On the exponential solution of differential equations for a linear operator}, Comm. Pure Appl. Math. \textbf{7}:649-673, 1954.

\bibitem{Mahony2}
R. Mahony, T. Hamel, and J.M. Pflimlin,  {\it Complementary Kalman filtering on matrix Lie groups}, IEEE Trans. Aut. Control {\bf 53}:1203-1218, 2008.

\bibitem{Mahony3}
M. Mallick, G.S. Chirikjian, and R. Mahony, {\it  A new numerical integration algorithm for rigid body dynamics on Lie groups}, IEEE Transactions on Robotics {\bf 30}:799-812, 2014.

\bibitem{MMM} 
J.C. Marrero, D. Martín de Diego and E. Mart\'inez, \textit{Discrete Lagrangian and Hamiltonian mechanics on Lie groupoids}, Nonlinearity \textbf{19}:1313-1348, 2006.


\bibitem{MaWe} 
J.E. Marsden and M. West, \textit{Discrete mechanics and variational integrators}, Acta Num. \textbf{10}:357-514, 2001.


\bibitem{McLac} R. McLachlan and G.R.W. Quispel, \textit{Splitting methods}, Acta Num. \textbf{11}:341-434, 2002.

\bibitem{McLa}
R.I. McLachlan, {\it Higher Order Symplectic Integrators},  Nonlinearity {\bf 8}:513-525, 1995.

\bibitem{munthe98} H. Munthe-Kaas, \textit{Runge-Kutta methods on Lie groups}, BIT Numerical Mathematics \textbf{38}:92-111, 1998.

\bibitem{munthe99} 
H. Munthe-Kaas, \textit{High order Runge-Kutta methods on manifolds}, J. Appl. Num. Maths. \textbf{29}:115-127, 1999.

\bibitem{munthe96}
H.Z. Munthe-Kaas and B. Owren, {\it On the construction of numerical methods preserving adiabatic invariants and related quantities},  Math. Comp. {\bf 65}:1353-1373, 1996.

\bibitem{Murray1}
R.M. Murray and S. Arimoto, {\it Numerical Solution of the Nonlinear Optimal Control Problem for Rigid Body Motion}, in {\sl Proceedings of the 29th IEEE Conference on Decision and Control 1990}, IEEE, Honolulu, 2000, pp. 2038-2043. 



\bibitem{odzijewicz00} 
A. Odzijewicz and A.M. Grundland, \textit{The Superposition Principle for the Lie Type first-order PDEs}, Rep. Math. Phys. \textbf{45}:293-306, 2000.


\bibitem{Pa_57}
 R.S. Palais, {\sl A Global Formulation of the Lie Theory of Transformation Groups}, 
 { Mem. Amer. Math. Soc.} {\bf 22}, Amer. Math. Soc., Providence, 1957.


\bibitem{PGG17}
A. Pandey, A. Ghose-Choudhury and P. Guha,
{\it 
Chiellini integrability and quadratically damped oscillators},
Int. J. Non-Linear Mech. {\bf 92}:153-159, 2017.

\bibitem{PW04}
A.V. Penskoi and P. Winternitz, {\it Discrete matrix Riccati equations with super- position formulas}, J. Math. Anal. Appl. {\bf 294}:533–547, 2004.


\bibitem{piet12} 
G. Pietrzkowski, 
\textit{Explicit solutions of the $ \mathfrak{a}_1 $-type Lie-Scheffers system and a general Riccati equation},  J. Dyn. Control Systems \textbf{18}:551-571, 2012.

\bibitem{Ros}
B.A. Rozenfel'd, 
{\sl A History of Non-Euclidean Geometry}, Springer, New York, 1988.

\bibitem{quarteroni07} 
A. Quarteroni, R. Sacco and F. Saleri, {\sl Numerical Mathematics}, Springer-Verlag, New York, 2007.

\bibitem{RW84}
D.W. Rand and P. Winternitz, 
{\it Nonlinear superposition principles: a new numerical method for solving matrix Riccati equations}, Comput. Phys. Comm. {\bf 33}:305-328, 1984.


\bibitem{Re72} 
W.T. Reid, {\sl Riccati Differential Equations}, Academic, New York, 1972.

\bibitem{SS} 
J.M. Sanz-Serna, 
\textit{Symplectic integrators for Hamiltonian problems: an overview}, Acta Num. {\bf 1}:{243-286}, 1992.

\bibitem{sardon15} 
C. Sardón, {\sl Lie systems, Lie symmetries and reciprocal transformations}, PhD Thesis, Universidad de Salamanca, Salamanca, 2015.

\bibitem{Murray2}
S. Sastry and R.M. Murray, {\it  Nonlinear feedback control for the stabilization of the attitude of a rigid body}, IEEE Trans. Aut. Control {\bf 32}:412-417, 1987.

\bibitem{Sattinger}
D.H. Sattinger and O.L. Weaver, {\sl Lie Groups and Algebras with Applications to Physics, Geometry, and Mechanics} Appl. Math. Sci. {\bf 61}, Springer-Verlag, Berlin, 1986. 

\bibitem{Sola}
J. Sola, J. Kuffner, and S.K. Agrawal, {\it A Quaternion-based Interpolation Technique for Smooth and Visually Appealing Motion Planning on Lie Groups} in: {\sl Proceedings of the 2018 IEEE International Conference on Robotics and Automation (ICRA)}, Brisbane, 2018, pp. 2181-2187.

\bibitem{sontag98} 
E.D. Sontag, {\sl Mathematical Control Theory: Deterministic Finite Dimensional Systems}, Springer-Verlag, New York, 1998.

\bibitem{varadajan}
V.S. Varadarajan, {\sl Lie Groups, Lie Algebras, and Their Representations}, Graduate Texts in Mathematics {\bf 102}, Springer, Los Angeles, 1985.

\bibitem{Verlet}
L. Verlet, {\it Computer "Experiments" on Classical Fluids. I. Thermodynamical Properties of Lennard-Jones Molecules}, Phys. Rev. {\bf 159}:98-103 (1967).

 \bibitem{Wi_83}
 P. Winternitz, {\it Lie groups and solutions of nonlinear differential equations} in {\sl Nonlinear phenomena(Oaxtepec, 1982)}, {Lect. Notes Phys.} {\bf 189}, Springer, Berlin, 1983, pp. 263-331.

 \bibitem{win82} 
P. Winternitz, \textit{Nonlinear action of Lie groups and superposition rules for nonlinear differential equations}, Phys. A {\bf 114}:105--113 (1982).
 
 \bibitem{Yaglom}
I.M. Yaglom, 
{\sl A simple non-Euclidean geometry and its physical basis}, Springer, New York, 1979.

\bibitem{Yoshida}
H. Yoshida, 
{\it Construction of Higher Order Symplectic Integrators}, Phys. Lett. A {\bf 150}:262-268, 1990.

\bibitem{Zanna} 
A. Zanna, \textit{Collocation and relaxed collocation for the Fer and Magnus expansions}, J. Numer. Anal. \textbf{36}:1145-1182, 1999.




\end{thebibliography}
\end{document}